\documentclass[reqno,12pt]{amsart}
\usepackage{amsmath,amssymb,latexsym,graphicx}
\def\1a{{\mathbf 1}}
\def\2a{{\mathbf 2}}
\def\3a{{\mathbf 3}}
\def\4a{{\mathbf 4}}
\def\5a{{\mathbf 5}}
\def\6a{{\mathbf 6}}
\def\7a{{\mathbf 7}}
\def\8a{{\mathbf 8}}
\def\9a{{\mathbf 9}}
\def\0a{\mathbf{10}}

\AtBeginDocument{

}

\numberwithin{equation}{section}
\newtheorem{theorem}{Theorem}[section]

\newtheorem{corollary}[theorem]{Corollary}
\newtheorem{definition}[theorem]{Definition}

\newtheorem{lemma}[theorem]{Lemma}

\newtheorem{example}[theorem]{Example}

\newcommand{\C}{\mathcal{C}}
\newcommand{\CI}{\mathcal{CI}}
\newcommand{\D}{\mathcal{D}}
\newcommand{\cP}{\mathcal{P}}
\newcommand{\ee}{\mathsf{e}}
\newcommand{\dd}{\mathsf{d}}

\newcommand{\nn}{\mathsf{n}}
\newcommand{\rc}{\mathrm{rc}}
\newcommand{\rev}{\mathrm{rev}}

\newcommand\GETOUT[1]{}
\newcommand{\s}{\mathfrak{S}}

\title[Centrosymmetric permutations]{Centrosymmetric Permutations and Involutions Avoiding $1243$ and $2143$}
\author{Mark F. Flanagan and Matteo Silimbani}
\address{School of Electrical, Electronic and Mechanical Engineering, University College Dublin, Ireland}
\email{mark.flanagan@ieee.org}
\address{Department of Mathematics, University of Bologna, Italy}
\email{silimban@dm.unibo.it}
\thanks{The first author was supported by Science Foundation Ireland (Grant no. 07/SK/I1252b).}
\begin{document}
\maketitle

\vspace{-5mm}
\begin{abstract}
A centrosymmetric permutation is one which is invariant under the reverse-complement operation, or equivalently one whose associated standard Young tableaux under the Robinson-Schensted algorithm are both invariant under the Sch\"{u}tzenberger involution. In this paper, we characterize the set of permutations avoiding $1243$ and $2143$ whose images under the reverse-complement mapping also avoid these patterns. We also characterize in a simple manner the corresponding Schr\"{o}der paths under a bijection of Egge and Mansour. We then use these results to enumerate centrosymmetric permutations avoiding the patterns $1243$ and $2143$. In a similar manner, centrosymmetric involutions avoiding these same patterns are shown to be enumerated by the Pell numbers.
\end{abstract}

\section{Introduction}

Let $\s_n$ denote the set of permutations of $\{ 1,2,\ldots, n\}$, and let $\s_n(\tau)$ denote the set of permutations of $\{ 1,2,\ldots, n\}$ which avoid the pattern $\tau$, i.e., which do not contain a subsequence order-isomorphic to $\tau$. More generally, let $\s_n(\tau_1,\tau_2,\ldots, \tau_k)$ denote the set of permutations of $\{ 1,2,\ldots, n\}$ which avoid all patterns $\tau_i$ for $i=1,2,\ldots, k$. 

For any permutation $\pi \in \s_n$, the \emph{reverse-complement} of $\pi$ is $\rc(\pi)=\pi'$ where $\pi'(i) = (n+1) - \pi(n+1-i)$ for each $i=1,2,\ldots, n$. Also, a permutation $\pi \in \s_n$ is said to be \emph{centrosymmetric} if and only if $\rc(\pi) = \pi$. Denote the set of centrosymmetric permutations by $\C_{n}$, and the set of centrosymmetric involutions by $\CI_{n}$. Also $\C_{n}(\tau)$, $\CI_{n}(\tau_1,\tau_2,\ldots, \tau_k)$ etc. are defined in the usual way.

It is well known that a permutation is centrosymmetric if and only if both of the standard Young tableaux yielded by the Robinson-Schensted algorithm are invariant under the Sch\"{u}tzenberger involution (see \cite{Schutzenberger} and \cite{Knuth} for more details). In work by Egge \cite{egge:restricted_symmetric_perms_L3}, permutations and involutions were enumerated which are invariant under the natural action of a subgroup of the symmetry group of a square. This included  enumeration and Wilf-equivalence classification of centrosymmetric permutations and involutions avoiding all patterns of length $3$. In a recent contribution also by the same author, $\left| \C_{2n}(k \; k-1 \cdots 2 1) \right|$ and $\left| \CI_{2n}(k \; k-1 \cdots 2 1) \right|$ were evaluated by counting self-evacuating standard Young tableaux and using the Robinson-Schensted correspondence \cite{egge:restricted_symmetric_perms_L4}. Other results along this line which have been achieved are the enumeration of the \emph{vexillary involutions}  (i.e., the set $\CI_{2n}(2143)$) by Guibert and Pergola \cite{Guibert}, and the set $\C_{2n}(123,2413)$ by Ostroff and Lonoff \cite{Lonoff_ostroff}. In \cite{Barnabei_bonetti_silimbani}, Barnabei \emph{et al.} enumerated many classes of pattern-avoiding centrosymmetric involutions by using a bijection with labeled Motzkin paths.

In this paper, centrosymmetric permutations and involutions avoiding the patterns $1243$ and $2143$ are enumerated. We begin by characterizing the set of permutations avoiding $1243$ and $2143$ whose images under the reverse-complement operation also avoid these patterns; to this end we make use of a result by Egge and Mansour \cite{egge_mansour} which puts permutations avoiding $1243$ and $2143$ in bijective correspondence with the set of Schr\"{o}der paths of an appropriate length. The characterization we require is particularly simple in the Schr\"{o}der path domain. We then enumerate centrosymmetric permutations avoiding $1243$ and $2143$ by enumerating the corresponding Schr\"{o}der paths. The corresponding enumeration for involutions is subsequently achieved by first proving that a Schr\"{o}der path $p$ corresponds to an involution under the bijection of Egge and Mansour if and only if $p$ is symmetric with respect to path reversal. In particular, the centrosymmetric involutions which avoid $1243$ and $2143$ are shown to be enumerated by the Pell numbers. 

It was proved in \cite{Albert} that the cardinality of the full class of permutations avoiding $1243$ and $2143$ whose reverse complements also avoid these patterns (i.e., the sequence $\left| \s_n(1243, 2143, 2134) \right|$) has a \emph{rational} generating function; it would be interesting to investigate whether similar methods can prove the rationality or otherwise of the generating functions of $\left| \C_{n}(1243,2143) \right|$ and $\left| \CI_{n}(1243,2143) \right|$.

\section{Permutations avoiding $1243$ and $2143$, and Schr\"{o}der paths}

The large Schr\"oder numbers $r_{n}$ are defined by $r_0=1$ and for all $n \ge 1$,
\begin{equation*}
r_{n} = r_{n-1} + \sum_{k=1}^{n} r_{k-1} r_{n-k}.
\end{equation*}
A \emph{Schr\"oder prefix} is a lattice path beginning at the point $(0,0)$ which may take only a finite number of steps from the set $\{\ee =(1,0), \nn =(0,1), \dd =(1,1) \}$ and which does not pass below the line $y=x$. Denote by $\mathcal{S}$ the set of all Schr\"oder prefixes. For $n \ge 1$, a \emph{Schr\"oder path} of length $n$ is a Schr\"oder prefix which terminates at the point $(n,n)$. Let $\mathcal{S}_n$ denote the set of Schr\"{o}der paths of length $n$. In this paper, such Schr\"{o}der paths will sometimes be denoted by the corresponding sequence of letters from $\{\ee, \nn, \dd\}$. Also denote by $\mathcal{S}_0$ the set containing the null path $\emptyset$ having length $0$. For any pair of Schr\"{o}der paths $p \in \mathcal{S}_m$ and $q \in \mathcal{S}_n$, denote by $p \, q \in \mathcal{S}_{m+n}$ the \emph{concatenation} of these Schr\"{o}der paths. The set $\mathcal{S}_n$ is enumerated by $r_n$ for $n \ge 0$. Also, the permutations $\s_{n+1}(1243,2143)$ for $n \ge 0$ are called \emph{Schr\"oder permutations} since they are enumerated by the large Schr\"oder numbers $r_n$. 
\begin{lemma}
A permutation $\pi$ lies in $\s_n(1243,2143)$ if and only if $\pi^{-1}$ lies in $\s_n(1243,2143)$.
\label{lemma:pi_inv}
\end{lemma}
\begin{proof}
This follows immediately from the observation that any occurrence of a pattern $\tau \in \{1243,2143\}$ in the permutation $\pi$ corresponds directly to an occurrence of the same pattern $\tau$ in $\pi^{-1}$.
\end{proof}

\begin{lemma}
Let $\pi_1, \pi_2 \in \s_n$. Then $\pi_1 = \rc(\pi_2)$ if and only if $\pi_1^{-1} = \rc(\pi_2^{-1})$.
\label{lemma:rc_inv}
\end{lemma}
\begin{proof}
This follows immediately from the observation that the substitution $j = \pi_1(i)$ allows the condition $\pi_1(i) + \pi_2(n+1-i) = n+1$ for all $i=1,2,\ldots,n$ to be rearranged as $\pi_1^{-1}(j) + \pi_2^{-1}(n+1-j)=n+1$ for all $j=1,2,\ldots,n$.
\end{proof}

\begin{definition}
For $t\in \{ 1,2,\ldots, n-1 \}$, a Schr\"oder path $p \in \mathcal{S}_n$ which contains an occurrence of $\dd$ joining $(t-1,t)$ to $(t,t+1)$ is said to have a \emph{level feature} at $t$, and a Schr\"oder path $p\in \mathcal{S}_n$ which contains an occurrence of $\ee\nn$ joining points $(t-1,t)$ to $(t,t)$ to $(t,t+1)$ is said to have a \emph{notch feature} at $t$. A Schr\"oder path $p\in \mathcal{S}_n$ is said to have a \emph{feature} at $t\in \{ 1,2,\ldots, n-1 \}$ if it has either a level feature or a notch feature at $t$. 
\label{def:features}
\end{definition}

\begin{definition}
For any Schr\"oder path $p\in \mathcal{S}_n$, the \emph{latest} (resp. level or notch) \emph{feature} of $p$ is the largest $t\in \{ 1,2,\ldots, n-1 \}$ at which $p$ contains a (resp. level or notch) feature, or is equal to $0$ if there is no such $t$. Also, for any Schr\"oder path $p\in \mathcal{S}_n$, the \emph{earliest} (resp. level or notch) \emph{feature} of $p$ is the smallest $t\in \{ 1,2,\ldots, n-1 \}$ at which $p$ contains a (resp. level or notch) feature, or is equal to $n$ if there is no such $t$. 
\label{def:latest_earliest_features}
\end{definition}

\begin{definition}
For $p \in \mathcal{S}_n$, the \emph{reversed path} $\rev(p)$ is the path obtained by applying the $\{ \ee,\nn,\dd \}$ steps of the path in reverse order, then replacing all occurrences of $\ee$ by $\nn$ and vice versa. Also for $p \in \mathcal{S}_n$, $\psi(p)$ is the path obtained by replacing all level features at $t\in \{ 1,2,\ldots, n-1 \}$ by notch features at $t$ and vice versa. 
\label{def:bijective_involutions}
\end{definition}

Egge and Mansour~\cite [\S 4]{egge_mansour} define a bijection $\varphi:\mathcal{S}_n \mapsto \s_{n+1}(1243,2143)$ from the set of Schr\"oder paths of length $n$ to the set of permutations {\bf of length} $\boldsymbol{n+1}$ avoiding $1243$ and $2143$. In the following we briefly describe the bijection $\varphi$; more details may be found in~\cite{egge_mansour}. An illustrative example of the bijection $\varphi$ will follow this formal description.

Let $p\in \mathcal{S}_n$ and let $s_i$ denote the transposition $(i,i+1)$ for each $i=1,2,\ldots, n$. We also use the convention $s_i s_j \pi = s_i(s_j(\pi))$.

\begin{description} 
\item[Step 1] Let $\Gamma_{r,s}$ denote the unit square which has diagonally opposite corners at $(r-1,s-1)$ and $(r,s)$. For each such square $\Gamma_{r,s}$ whose top-left corner lies below the path $p$ and above the line $y=x$, place the label $r$ in the top-left corner. We will next construct a sequence of permutations $\sigma_i$, $i=1,2,\ldots,k$; initially set $i=1$.

\item[Step 2] Locate the labeled square $\Gamma_{r,s}$ with largest label $r$. The permutation $\sigma_i$ is equal to the sequence of transpositions $s_t$ where the subscript $t$ sequentially takes on all label values starting with this $r$ and continuing diagonally to the lower left until an $\nn$ step is encountered. Remove the labels which were used as subscripts. If there are now no more labeled squares, we are finished; otherwise set $i$ to $i+1$ and repeat Step 2.

\item[Step 3]
The image $\varphi(p)$ of the Schr\"oder path $p$ under the bijection $\varphi$ is defined as
\begin{equation}
\varphi(p) = \sigma_k \sigma_{k-1} \cdots \sigma_1 (n+1,n,n-1,\ldots, 2,1) \; .
\label{eq:Egge_Mansour_bijection}
\end{equation}
\end{description}

\begin{figure}
\begin{center}\includegraphics[%
  width=0.5\columnwidth,
  keepaspectratio]{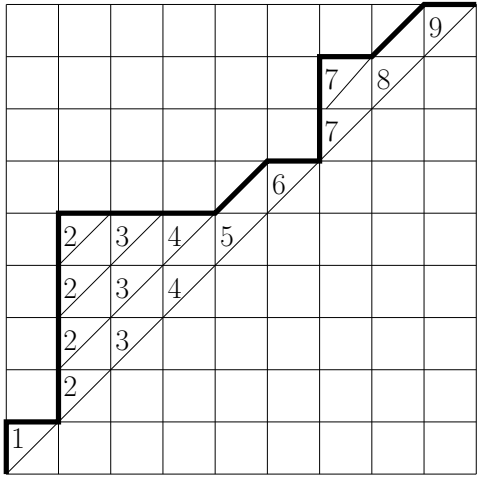}
\end{center}
\caption{\label{cap:Schroder1}The Schr\"oder path $p$ of Example \ref{ex:running_example}.}
\end{figure}

The following example illustrates Definitions \ref{def:features} and \ref{def:latest_earliest_features} as well as the bijection $\varphi$, and will also serve as a useful example when we consider the reverse-complements of permutations in $\s_{n+1}(1243,2143)$.

\begin{example}\label{ex:running_example}
Consider the path $p = \nn\ee\nn\nn\nn\nn\ee\ee\ee\dd\ee\nn\nn\ee\dd\ee \in \mathcal{S}_9$ shown in Figure \ref{cap:Schroder1}.
The path $p$ has an (earliest) level feature at $t=5$ and another (latest) level feature at $t=8$. Also, the path $p$ has an (earliest) notch feature at $t=1$ and another (latest) notch feature at $t=6$. To compute $\varphi(p)$, note that here we have
$\sigma_1=s_9s_8s_7$, $\sigma_2=s_7$, $\sigma_3= s_6s_5s_4s_3s_2$, $\sigma_4=s_4s_3s_2$, $\sigma_5=s_3s_2$, $\sigma_6=s_2$ and $\sigma_7=s_1$. So
\begin{eqnarray*}
\varphi(p) &=& \sigma_7 \sigma_6 \sigma_5 \sigma_4 \sigma_3 \sigma_2 \sigma_1 (10,9,8,7,6,5,4,3,2,1) \\
&=& s_1 \; s_2 \; s_3s_2 \; s_4s_3s_2 \; s_6s_5s_4s_3s_2 \; s_7 \; s_9s_8s_7 (10,9,8,7,6,5,4,3,2,1) \\
&=& (5,10,6,7,8,2,9,3,1,4) \in \s_{10}(1243,2143).
\end{eqnarray*}
A detailed evolution of $(10,9,8,7,6,5,4,3,2,1)$ to $\varphi(p)$ under the application of the permutations $\sigma_i$ for $i=1,2,\ldots,k$ is provided in Table \ref{ttable_1}. 
\end{example}

\begin{lemma}
A Schr\"oder path $p \in \mathcal{S}_n$ contains an occurrence of $\dd$ joining $(t-1,t-1)$ to $(t,t)$ for $t \in \{ 1,2,\ldots,n \}$ if and only if the largest $t$ numbers in $\{ 1,2,\ldots, n+1\}$ occupy the first $t$ positions of $\pi = \varphi(p)$.
\label{lemma:no_mixing}
\end{lemma}
\begin{proof}
With reference to (\ref{eq:Egge_Mansour_bijection}), note that $\pi = \varphi(p)$ is obtained by applying the permutations $\sigma_1, \sigma_2, \ldots, \sigma_k$ to a permutation in which the largest $t$ numbers in $\{ 1,2,\ldots, n+1\}$ occupy the first $t$ positions. Next observe that the path $p \in \mathcal{S}_n$ contains an occurrence of $\dd$ joining $(t-1,t-1)$ to $(t,t)$ if and only there is no occurrence of $s_t$ among the permutations $\sigma_1, \sigma_2, \ldots, \sigma_k$. This is exactly the condition under which there will be no ``mixing'' of the two parts of the permutation, i.e., the largest $t$ numbers in $\{ 1,2,\ldots, n+1\}$ will remain in the first $t$ positions of $\pi = \varphi(p)$. 
\end{proof}

\begin{theorem}
If $p$ is a Schr\"{o}der path of length $n \ge 1$, and $\pi = \varphi(p)$, then $\varphi(\rev(p)) = \pi^{-1}$.
\label{lemma:inverse_through_EM_bijection}
\end{theorem}
\begin{proof}
We prove this by induction on the length $n$ of the path $p$. The base case $n=1$ holds trivially. Next, assume the result holds for all paths of length at most $n-1$ (where $n > 1$), and consider a path $p$ of length $n$. Denote $\pi = \varphi(p)$ and $\tilde{\pi} = \varphi(\rev(p))$. Our goal is to show that $\tilde{\pi} = \pi^{-1}$.

If no point $(t,t)$ lies on the path $p$ for $t = 1,2,\ldots,n-1$, then we may write $p = \nn \, p_1 \, \ee$ for some Schr\"{o}der path $p_1$ of length $n-1$. Let $\pi_1 = \varphi(p_1)$. Then, with reference to obtaining $\pi = \varphi(p)$ from $(n+1,n,n-1,\ldots,2,1)$ via (\ref{eq:Egge_Mansour_bijection}), first $\sigma_1 = s_n s_{n-1} \cdots s_2 s_1$ is applied; this is a cyclic shift which yields $\pi(n+1) = n+1$. The subsequent sequence of permutations $\{\sigma_i \; : \; i \ge 2 \}$ is the same as that involved in obtaining $\varphi(p_1)$ from $(n,n-1,\ldots,2,1)$ via (\ref{eq:Egge_Mansour_bijection}); thus $\pi(i) = \pi_1(i)$ for all $i=1,2,\ldots,n$. Since $\rev(p) = \nn \, \rev(p_1) \, \ee$, using the same argument together with the induction hypothesis yields $\tilde{\pi}(i) = \pi^{-1}_1(i)$ for all $i=1,2,\ldots,n$, and $\tilde{\pi}(n+1) = n+1$. Thus $\tilde{\pi} = \pi^{-1}$. 

\begin{table}
\begin{center}
\caption{\label{ttable_1}
Illustration of the bijection $\varphi(p)$ for the Schr\"oder path $p \in \mathcal{S}_9$ given in Example \ref{ex:running_example}. The table shows the evolution of the permutation from $(10,9,8,7,6,5,4,3,2,1)$ towards $\varphi(p)$ as each permutation $\sigma_i$ is applied for $i=1,2,\ldots,k$.}

{\small
 $\begin{array}{|c|c|}
 \hline
\begin{array}[t]{c}
\mathrm{Permutation} \\
\hline 
\mathrm{Start} \\
\sigma_1 = s_9 s_8 s_7 \\
\sigma_2 = s_7 \\
\sigma_3 = s_6 s_5 s_4 s_3 s_2 \\
\sigma_4 = s_4 s_3 s_2 \\
\sigma_5 = s_3 s_2 \\
\sigma_6 = s_2 \\
\sigma_7= s_1 \\
\end{array} &
\begin{array}[t]{c}
\mathrm{Result} \\
\hline 
10 \; 9 \; 8 \; 7 \; 6 \; 5 \; 4 \; 3 \; 2 \; 1 \\
10 \; 9 \; 8 \; 7 \; 6 \; 5 \; \3a \; \2a \; \1a \; \4a \\
10 \; 9 \; 8 \; 7 \; 6 \; 5 \; \2a \; \3a \; 1 \; 4 \\
10 \; \8a \; \7a \; \6a \; \5a \; \2a \; \9a \; 3 \; 1 \; 4 \\
10 \; \7a \; \6a \; \5a \; \8a \; 2 \; 9 \; 3 \; 1 \; 4 \\
10 \; \6a \; \5a \; \7a \; 8 \; 2 \; 9 \; 3 \; 1 \; 4 \\
10 \; \5a \; \6a \; 7 \; 8 \; 2 \; 9 \; 3 \; 1 \; 4 \\
\5a \; \0a \; 6 \; 7 \; 8 \; 2 \; 9 \; 3 \; 1 \; 4 \\
\end{array} \\
\hline
\end{array}$
}
\end{center}
\end{table}

If $p$ terminates in a $\dd$ step, then we may write $p = p_1 \, \dd$ for some Schr\"{o}der path $p_1$ of length $n-1$. Let $\pi_1 = \varphi(p_1)$. Then, the sequence of permutations $\{\sigma_i\}$ involved in obtaining $\varphi(p)$ from $(n+1,n,n-1,\ldots,2,1)$ via (\ref{eq:Egge_Mansour_bijection}) is the same as that involved in obtaining $\varphi(p_1)$ from $(n,n-1,\ldots,2,1)$ via (\ref{eq:Egge_Mansour_bijection}); thus $\pi(i) = \pi_1(i) + 1$ for all $i=1,2,\ldots,n$, and $\pi(n+1) = 1$. Also, since $\rev(p) = \dd \, \rev(p_1)$, the sequence of permutations used to obtain $\varphi(\rev(p))$ from $(n+1,n,n-1,\ldots,2,1)$ via (\ref{eq:Egge_Mansour_bijection}) is the same as that used to obtain $\varphi(\rev(p_1))$ from $(n,n-1,\ldots,2,1)$ via (\ref{eq:Egge_Mansour_bijection}), but with each $s_t$ replaced by $s_{t+1}$. Together with the induction hypothesis, this yields $\tilde{\pi}(i) = \pi^{-1}_1(i-1)$ for all $i=2,3,\ldots,n+1$, and $\tilde{\pi}(1) = n+1$. Thus $\tilde{\pi} = \pi^{-1}$.
 
Therefore, hereafter we need only prove the inductive step for Schr\"{o}der paths $p$ of length $n$ which terminate in an $\ee$ step and such that the point $(t,t)$ lies on $p$ for some $t \in \{ 1,2,\ldots,n-1 \}$. For such paths, there exists a largest value of $t \in \{1,2,\ldots,n-1\}$ for which $(t,t)$ lies on the path $p$; denote by $n-i$ this largest value, where $i \in \{ 1,2,\ldots,n-1 \}$. Note that an $\nn$ step must connect $(n-i,n-i)$ to $(n-i,n-i+1)$. 

First consider the case where the point $(n-i,n-i)$ on $p$ is reached via an $\ee$ step. Note that we may write $p = q \, r$ where $q$ and $r$ are Schr\"oder paths of length $n-i$ and $i$ respectively. 

With reference to (\ref{eq:Egge_Mansour_bijection}), the permutation $\pi_2=\varphi(q)$ is obtained by applying some sequence of permutations $\{ \tau_j \}$ to the permutation $(n+1-i,n-i,n-1-i,\ldots, 2,1)$. Similarly, the permutation $\pi_1=\varphi(r)$  is obtained by applying some sequence of permutations $\{ \sigma_j \}$ to the permutation $(i+1,i,i-1,\ldots, 2,1)$. The image of $p$ under the bijection $\varphi$ is then obtained by applying the sequence of permutations $\{ \nu_j \}$, followed by the sequence of permutations $\{ \tau_j \}$, to the permutation $(n+1,n,n-1,\ldots, 2,1)$, where the sequence $\{ \nu_j \}$ is simply the sequence $\{ \sigma_j \}$ with each transposition $s_t$ replaced by $s_{t+n-i}$. 

The application of the permutations $\{ \nu_j\}$ (corresponding to the path $r$) results in the final $i+1$ values of $(n+1,n,n-1,\ldots, 2,1)$ being replaced with the permutation $\pi_1$. Then, the application of the permutations $\{ \tau_j \}$ (corresponding to the path $q$) results in the initial $n+1-i$ values of the resulting permutation being replaced with $\pi_2 + i$, except for the unique position $t\in \{1,2,\ldots,n+1-i\}$ with $\pi_2(t) = 1$, which takes the value $\pi_1(1)$.

We may summarize the results of the previous paragraph as follows: for all $t\in\{1,2,\ldots,n+1\}$, 
\begin{equation}
\pi(t) = \left\{ \begin{array}{ll}
\pi_1(1) & \textrm{ if } \pi_2(t) = 1 \\
\pi_1(t-n+i) & \textrm{ if } t > n+1-i \\
\pi_2(t) + i & \textrm{ otherwise. }\end{array}\right. \;
\label{eq:case_1_pi}
\end{equation}

Since $\rev(p) = \rev(r) \, \rev(q)$, applying the same argument and invoking the induction hypothesis yields that, for all $s\in\{1,2,\ldots,n+1\}$, 
\begin{equation}
\tilde{\pi}(s) = \left\{ \begin{array}{ll}
\pi_2^{-1}(1) & \textrm{ if } \pi_1^{-1}(s) = 1 \\
\pi_2^{-1}(s-i) & \textrm{ if } s > i+1 \\
\pi_1^{-1}(s) + n-i & \textrm{ otherwise. }\end{array}\right. \;
\label{eq:case_1_pi_tilde}
\end{equation}
The reader may check that \eqref{eq:case_1_pi} and \eqref{eq:case_1_pi_tilde} together imply (through the identification $\pi(t) = s$) that $\tilde{\pi} = \pi^{-1}$, as required.

Finally, consider the case where the point $(n-i,n-i)$ on $p$ is reached via a $\dd$ step. In this case, we may write $p = q \, \dd \, r$, where $q$ and $r$ are Schr\"oder paths of length $n-i-1$ and $i$ respectively. Again denote $\pi_1=\varphi(r)$ and $\pi_2=\varphi(q)$. Arguing in a similar manner to the previous case, we find that in obtaining $\pi=\varphi(p)$ via (\ref{eq:Egge_Mansour_bijection}), the sequence of permutations corresponding to the path $r$ results in the final $i+1$ values of $(n+1,n,n-1,\ldots, 2,1)$ being replaced with the permutation $\pi_1$. Then, the sequence of permutations corresponding to the path $q$ results in the initial $n-i$ values of the resulting permutation being replaced with $\pi_2 + i + 1$. Summarizing,
\begin{equation}
\pi(t) = \left\{ \begin{array}{ll}
\pi_1(t-n+i) & \textrm{ if } t > n-i \\
\pi_2(t) + i + 1 & \textrm{ otherwise. }\end{array}\right. \;
\label{eq:case_2_pi}
\end{equation}
Since $\rev(p) = \rev(r) \, \dd \, \rev(q)$, applying the same argument and invoking the induction hypothesis yields that, for all $s\in\{1,2,\ldots,n+1\}$, 
\begin{equation}
\tilde{\pi}(s) = \left\{ \begin{array}{ll}
\pi_2^{-1}(s-i-1) & \textrm{ if } s > i + 1 \\
\pi_1^{-1}(s) + n-i & \textrm{ otherwise. }\end{array}\right. \;
\label{eq:case_2_pi_tilde}
\end{equation}
As before, \eqref{eq:case_2_pi} and \eqref{eq:case_2_pi_tilde} together imply (through the identification $\pi(t) = s$) that $\tilde{\pi} = \pi^{-1}$, as required. The result then follows by the principle of induction.
\end{proof}

\begin{corollary}\label{cor:involutions}
Let $p \in \mathcal{S}_n$ and $\pi = \varphi(p)$. Then $\pi$ is an involution if and only if $\rev(p) = p$.
\end{corollary}
It follows that the number of involutions on $\{1,2,3,\ldots,n+1\}$ which avoid $1243$ and $2143$ is equal to the number of Schr\"{o}der paths $p \in \mathcal{S}_n$ which are symmetric with respect to path reversal. This latter fact may also be deduced from the bijection given in \cite[\S 2]{deng_dukes:Schroder_involutions}. 

The following example illustrates the inductive step in the proof of Theorem \ref{lemma:inverse_through_EM_bijection}.
\begin{example}\label{ex:inverse_through_EM_bijection}
Consider the Schr\"{o}der path $p = \nn\dd\nn\dd\nn\ee\ee\ee\nn\nn\nn\ee\nn\ee\ee\ee \in \mathcal{S}_9$ illustrated in Figure \ref{cap:Schroder_general}. The path $p$ terminates in an $\ee$ step, and the largest $t \in \{ 1,2,\ldots,n-1 \} = \{ 1,2,\ldots,8 \}$ for which the point $(t,t)$ lies on $p$ is given by $t = n-i = 5$, i.e., $i=4$. The point $(5,5)$ on $p$ is reached via an $\ee$ step. Here we have $p = q \, r$ where $q = \nn\dd\nn\dd\nn\ee\ee\ee \in \mathcal{S}_5$ and $r = \nn\nn\nn\ee\nn\ee\ee\ee \in \mathcal{S}_4$. The reader may verify that $\pi_2 = \varphi(q) = (5,3,1,2,4,6)$ and that $\varphi(\rev(q)) = \varphi(\nn\nn\nn\ee\dd\ee\dd\ee) = (3,4,2,5,1,6) = \pi_2^{-1}$. Also $\pi_1 = \varphi(r) = (1,3,2,4,5)$ and since $\rev(r) = r$ we have $\varphi(\rev(r)) = \pi_1$; the reader may verify that $\pi_1 = \pi_1^{-1}$.

\begin{figure}
\begin{center}\includegraphics[%
  width=0.5\columnwidth,
  keepaspectratio]{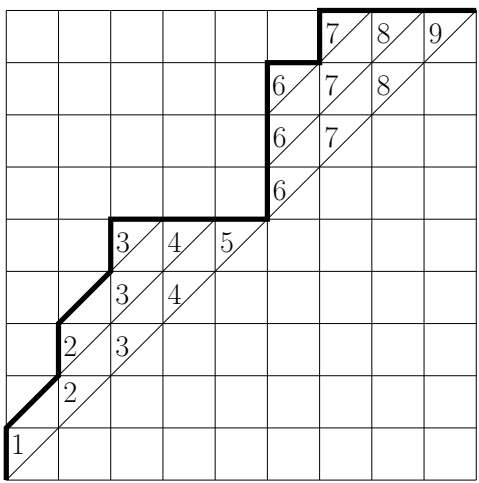}
\end{center}
\caption{\label{cap:Schroder_general}The Schr\"oder path $p$ of Example \ref{ex:inverse_through_EM_bijection}.}
\end{figure}

The application of the first set of permutations $\{ \nu_j\}$ (corresponding to the path $r$) results in the final $i+1=5$ values of $(10,9,8,7,6,5,4,3,2,1)$ being replaced with $\pi_1 = (1,3,2,4,5)$ (c.f. row 2 of Table \ref{ttable_general}). Then, the application of the second set of permutations $\{ \tau_j \}$ (corresponding to the path $q$) results in the initial $n+1-i=6$ values of the resulting permutation being replaced with $\pi_2 + i = \pi_2 + 4 = (9,7,5,6,8,10)$, except for position $\pi_2^{-1}(1) = 3$ which takes the value $\pi_1(1) = 1$ (c.f. row 3 of Table \ref{ttable_general}). Thus \eqref{eq:case_1_pi} holds in this case.
\end{example}
\begin{table}[h]
\begin{center}
\caption{\label{ttable_general}
Illustration of the steps involved in obtaining $\pi=\varphi(p)$ for the Schr\"oder path $p \in \mathcal{S}_9$ given in Example \ref{ex:inverse_through_EM_bijection}. The table shows the evolution of the permutation from the starting point of $(10,9,8,7,6,5,4,3,2,1)$.}
{\small
 $\begin{array}{|c|c|}
 \hline
\begin{array}[t]{c}
\mathrm{Permutation} \\
\hline 
\mathrm{Start} \\
\{ \nu_j \} \\
\{ \tau_j \} \\
\end{array} &
\begin{array}[t]{c}
\mathrm{Result} \\
\hline 
10 \; 9 \; 8 \; 7 \; 6 \; 5 \; 4 \; 3 \; 2 \; 1 \\
10 \; 9 \; 8 \; 7 \; 6 \; \underline{\1a} \; \3a \; \2a \; \4a \; \5a \\
\9a \; \7a \; \underline{\1a} \; \6a \; \8a \; \0a \; 3 \; 2 \; 4 \; 5 \\
\end{array} \\
\hline
\end{array}$
}
\end{center}
\end{table}
\section{Centrosymmetric permutations and involutions avoiding $1243$ and $2143$}
\begin{definition}
The set $\cP$ is the set of Schr\"oder paths formed by concatenating a finite sequence of elements from $\dd \cup \{ \nn^k \ee^k \mid k > 0 \}$. The set $\D_n \subseteq \mathcal{S}_n$ is the set of Schr\"oder paths of length $n$ formed by concatenating a finite sequence of elements from $\dd \cup \nn \cP \ee$.
\label{def:peaks_etc}
\end{definition}
Less formally, the set $\D_n$ is the set of Schr\"{o}der paths of length $n$ which, whenever they rise above the ``main superdiagonal'' $y = x + 1$, do so by a sequence of $\nn$ steps immediately followed by an equal number of $\ee$ steps. 
\begin{lemma}
Let $p \in \D_n$ be a Schr\"{o}der path with no features, and which does not contain a $\dd$ step connecting $(t-1,t-1)$ to $(t,t)$ for any $t \in \{1,2,\ldots,n\}$. Furthermore, let $\varphi(p) = \pi$, $p' = \rev (\psi(p))$, and $\varphi(p') = \pi'$. Then $\pi' = \rc(\pi)$, and $\pi(1) = 1$. 
\label{lemma:induction_base_case}
\end{lemma}
\begin{proof}
The Schr\"{o}der path $p$ is of the form $\nn \alpha_1 \alpha_2 \cdots \alpha_{r-2} \alpha_{r-1} \ee$ where $\alpha_j = \nn^{m_{j+1}-m_{j}} \ee^{m_{j+1}-m_{j}}$ for some $r$ numbers $1 = m_1 < m_2 < \cdots < m_{r-1} < m_r = n$. In terms of the permutations $\{\sigma_i\}$ involved in obtaining $\varphi(p)$ via (\ref{eq:Egge_Mansour_bijection}), note that $\alpha_j$ ($1 \le j < r$) corresponds to a \emph{reversal} of the numbers in positions $\{ m_j, m_j+1, \ldots, m_{j+1} \}$; denote this reversal by $\delta(m_j, m_{j+1})$. Thus
\begin{equation}
\varphi(p) = \pi = \delta(m_1, m_2) \delta(m_2, m_3) \cdots \delta(m_{r-1}, m_r) \sigma_1 (n+1,n,n-1,\ldots, 2,1) \; ,
\label{eq:base_case_step_1}
\end{equation}
where $\sigma_1 = s_n s_{n-1} \cdots s_2 s_1$.

Since $p$ has no features we have 
\[
p' = \rev (\psi(p)) = \rev (p) = \nn \alpha_{r-1} \alpha_{r-2} \cdots \alpha_2 \alpha_1 \ee \; ,
\] 
and thus
\begin{align}
&\varphi(p') = \pi' = \delta(n+1-m_r, n+1-m_{r-1}) \delta(n+1-m_{r-1}, n+1-m_{r-2}) \nonumber \\
&\cdots \delta(n+1-m_{2}, n+1-m_{1}) \sigma_1 (n+1,n,n-1,\ldots, 2,1) \; .
\label{eq:base_case_step_2}
\end{align}
From comparison of (\ref{eq:base_case_step_1}) and (\ref{eq:base_case_step_2}), it may be observed that $\pi' = \rc(\pi)$ due to the mirror symmetry of the applied reversal operations, together with the fact that $\pi(n+1) = \pi'(n+1) = n+1$ and $\pi(1) = \pi'(1) = 1$.
\end{proof}
\begin{figure}
\begin{center}
  \includegraphics[%
  width=0.48\columnwidth,
  keepaspectratio]{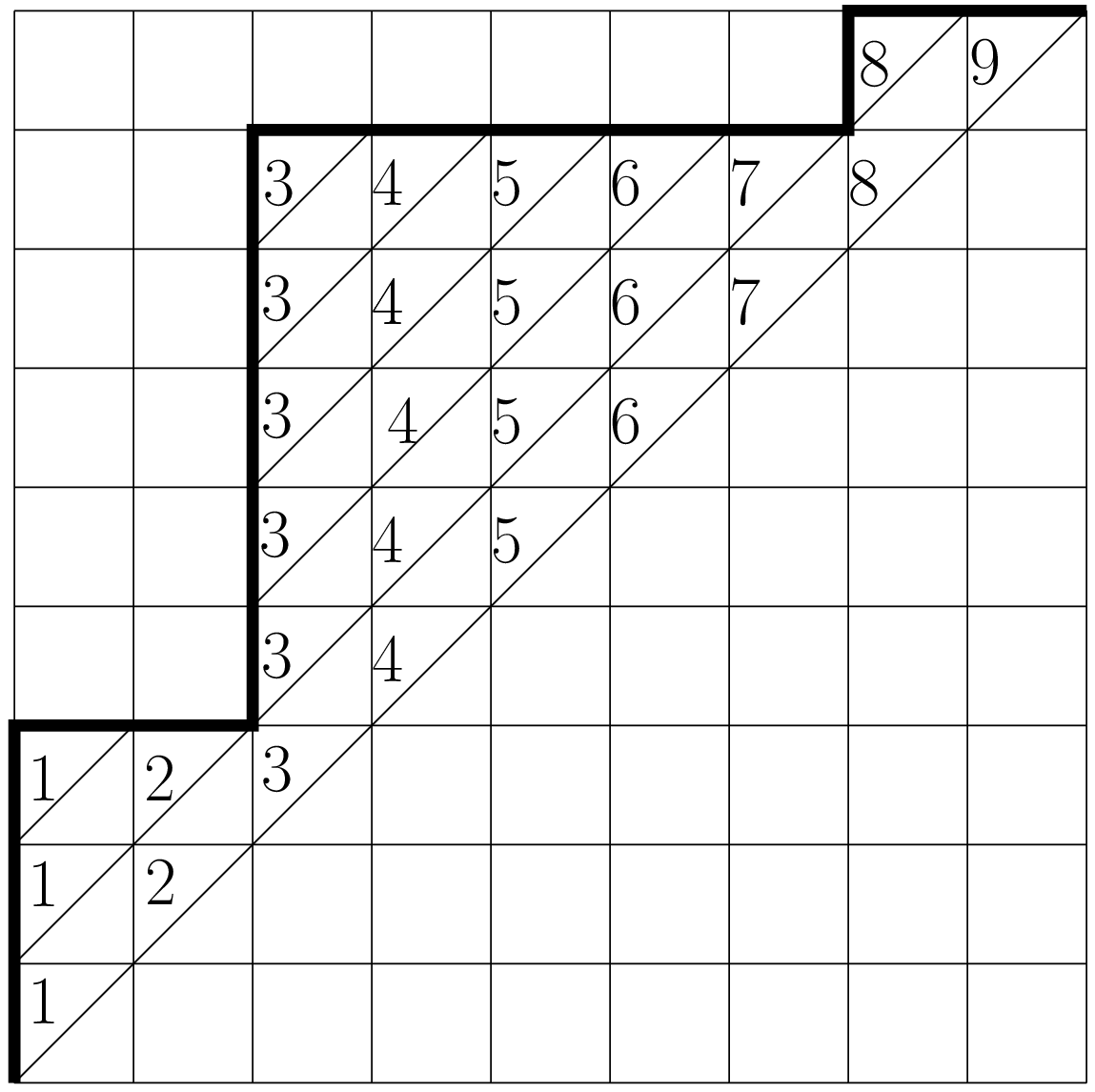}
  \includegraphics[%
  width=0.48\columnwidth,
  keepaspectratio]{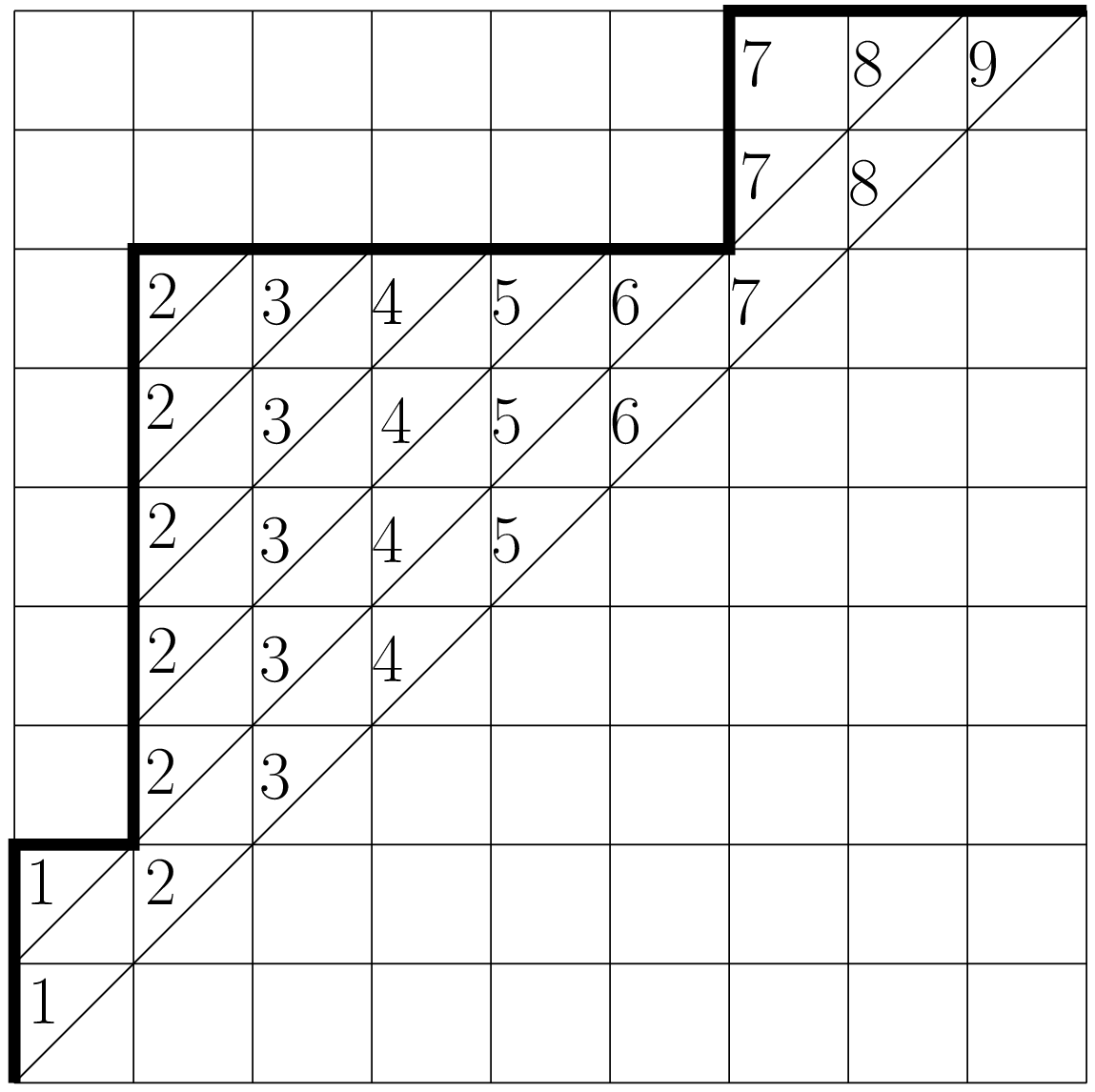}
\vspace{2mm} \\
(a) \hspace{5cm} (b) 
\end{center}
\caption{\label{cap:Schroder_no_features1} (a) The Schr\"oder path $p \in \D_9$ and (b) the Schr\"oder path $p' = \rev (\psi(p))$ of Example \ref{ex:base_case}.}
\end{figure}
\begin{example}
Consider the Schr\"{o}der path $p \in \D_9$ illustrated in Figure \ref{cap:Schroder_no_features1} (a). The path $p$ has no features, and does not contain a $\dd$ step connecting $(t-1,t-1)$ to $(t,t)$ for any $t \in \{1,2,\ldots,9\}$. The path $p' = \rev (\psi(p)) = \rev(p)$ is shown in Figure \ref{cap:Schroder_no_features1} (b). Table \ref{ttable_3} shows the sequence of reversals applied to the permutation $(10,9,8,7,6,5,4,3,2,1)$ to obtain $\pi = \varphi(p)$ (upper Table \ref{ttable_3}) and $\pi' = \varphi(p')$ (lower Table \ref{ttable_3}) respectively. It may be observed that $\pi' = \rc(\pi)$ holds due to the mirror symmetry of the applied reversal operations as well as the fact that both permutations $\pi$ and $\pi'$ interchange the values $1$ and $n+1=10$.
\label{ex:base_case}
\end{example}
\begin{lemma}\label{lemma:concatenation}
Let $p \in \D_n$ be a Schr\"{o}der path which does not contain a $\dd$ step connecting $(t-1,t-1)$ to $(t,t)$ for any $t \in \{1,2,\ldots,n\}$, and which contains no level features. If $\pi = \varphi(p)$, then $\pi(1) = 1$.
\end{lemma}
\begin{proof}
Write $p = q_1 q_2 \cdots q_l$ where each Schr\"{o}der path $q_i$ is of the form described in Lemma \ref{lemma:induction_base_case}; therefore, each permutation $\pi_i = \varphi(q_i)$ satisfies $\pi_i(1) = 1$. Letting $m_i>0$ denote the length of each path $q_i$, we have $\sum_{j=1}^{l} m_j = n$. Then, with reference to obtaining $\varphi(p)$ from $(n+1,n,n-1,\ldots,2,1)$ via \eqref{eq:Egge_Mansour_bijection}, permutations are applied corresponding to each path $q_i$ for $i=l,l-1,\ldots,2,1$ respectively, in each case moving the entry $1$ from position $\sum_{j=1}^i m_j + 1$ to position $\sum_{j=1}^{i-1} m_j + 1$. Therefore, after all permutations are applied we must have $\pi(1) = 1$.
\end{proof}
\begin{table}
\begin{center}
\caption{\label{ttable_3}
Illustration of the subsequence reversal operations involved in obtaining $\pi=\varphi(p)$ (upper) and $\pi' = \varphi(p')$ (lower) for the Schr\"oder paths $p,p' \in \D_9$ given in Example \ref{ex:base_case}. The table shows the evolution of the permutation from $(10,9,8,7,6,5,4,3,2,1)$ with the initial cyclic shift $\sigma_1$ followed by the reversals $\delta(\cdot,\cdot)$ in each case.}

{\small
 $\begin{array}{|c|c|}
 \hline
\begin{array}[t]{c}
\mathrm{Permutation} \\
\hline 
\mathrm{Start} \\
\sigma_1 = s_9 s_8 s_7 s_6 s_5 s_4 s_3 s_2 s_1 \\
\delta(8,9) \\
\delta(3,8) \\
\delta(1,3) \\
\end{array} &
\begin{array}[t]{c}
\mathrm{Result} \\
\hline 
10 \; 9 \; 8 \; 7 \; 6 \; 5 \; 4 \; 3 \; 2 \; 1 \\
\9a \; \8a \; \7a \; \6a \; \5a \; \4a \; \3a \; \2a \; \1a \; \0a \\
9 \; 8 \; 7 \; 6 \; 5 \; 4 \; 3 \; \1a \; \2a \; 10 \\
9 \; 8 \; \1a \; \3a \; \4a \; \5a \; \6a \; \7a \; 2 \; 10 \\
\1a \; \8a \; \9a \; 3 \; 4 \; 5 \; 6 \; 7 \; 2 \; 10 \\
\end{array} \\
\hline
\end{array}$
}
{\small
 $\begin{array}{|c|c|}
 \hline
\begin{array}[t]{c}
\mathrm{Permutation} \\
\hline 
\mathrm{Start} \\
\sigma_1 = s_9 s_8 s_7 s_6 s_5 s_4 s_3 s_2 s_1 \\
\delta(7,9) \\
\delta(2,7) \\
\delta(1,2) \\
\end{array} &
\begin{array}[t]{c}
\mathrm{Result} \\
\hline 
10 \; 9 \; 8 \; 7 \; 6 \; 5 \; 4 \; 3 \; 2 \; 1 \\
\9a \; \8a \; \7a \; \6a \; \5a \; \4a \; \3a \; \2a \; \1a \; \0a \\
9 \; 8 \; 7 \; 6 \; 5 \; 4 \; \1a \; \2a \; \3a \; 10 \\
9 \; \1a \; \4a \; \5a \; \6a \; \7a \; \8a \; 2 \; 3 \; 10 \\
\1a \; \9a \; 4 \; 5 \; 6 \; 7 \; 8 \; 2 \; 3 \; 10 \\
\end{array} \\
\hline
\end{array}$
}
\end{center}
\end{table}
\begin{lemma}\label{thm:earliest_level_feature}
Let $p \in \D_n$ be a Schr\"{o}der path which does not contain a $\dd$ step connecting $(t-1,t-1)$ to $(t,t)$ for any $t \in \{1,2,\ldots,n\}$. If the earliest level feature of $p$ is at $t=k$, then $\pi(1) = n-k+1$.
\end{lemma}
\begin{proof}

Let $q \in \mathcal{D}_k$ denote the Schr\"{o}der path of length $k$ obtained by following $p$ up the point $(k-1,k)$ and terminating via an $\ee$ step at the point $(k,k)$. Note that $q$ is of the form described in Lemma \ref{lemma:concatenation}, and thus $\pi_1 = \varphi(q)$ satisfies $\pi_1(1) = 1$. Next, with reference to obtaining $\varphi(p)$ from $(n+1,n,n-1,\ldots,2,1)$ via \eqref{eq:Egge_Mansour_bijection}, initially the entry $n-k+1$ lies in position $k+1$. First a sequence of permutations is applied which uses only transpositions $s_t$ with $t > k+1$; after this, the entry $n-k+1$ remains in position $k+1$. Then, a sequence of transpositions is applied which contains $s_k$; this moves the entry $n-k+1$ into position $k$. Immediately after this, a sequence of permutations is applied containing only transpositions $s_t$ with $t>k$; during this process, entry $n-k+1$ remains in position $k$. For the remaining sequence of permutations, observe that since entry $n-k+1$ starts in position $k$, its final position is the same as the final position of entry $1$ when obtaining $\pi_1 = \varphi(q)$ via \eqref{eq:Egge_Mansour_bijection}, i.e., $\pi(n-k+1) = 1$.
\end{proof}
\begin{example}
Consider the Schr\"{o}der path $p = \nn\nn\ee\nn\ee\ee\nn\nn\ee\dd\nn\ee\ee\nn\nn\ee\ee \in \D_9$ illustrated in Figure \ref{cap:Schroder_concat}. The earliest level feature of $p$ is at $t=k=5$. Here we have $q = \nn\nn\ee\nn\ee\ee\nn\nn\ee\ee$ and $\pi_1 = \varphi(q) = (1,5,4,6,2,3)$. The Schr\"{o}der path $q$ is of the form described in Lemma \ref{lemma:concatenation}, and accordingly $\pi_1(1) = 1$. Table \ref{ttable_concat} summarizes the steps involved in obtaining $\pi = \varphi(p)$ from $(10,9,8,7,6,5,4,3,2,1)$ via \eqref{eq:Egge_Mansour_bijection}. The permutations $\{ \sigma_i \}$ prior to that containing transposition $s_k = s_5$ involve only transpositions $s_8$ and $s_9$, and thus do not affect the position of entry $n-k+1 = 5$. Next, the permutation $s_7 s_6 s_5 s_4$ is applied which moves entry $n-k+1 = 5$ to position $k=5$. The next permutation, $s_6$, does not alter the position of the entry $5$. The final sequence of permutations is simply that used to obtain $\pi_1 = \varphi(q)$ via \eqref{eq:Egge_Mansour_bijection}, but with $\sigma_1 = s_5 s_4$ omitted; therefore we finally obtain $\pi(1) = 5$.
\begin{figure}
\begin{center}\includegraphics[%
  width=0.5\columnwidth,
  keepaspectratio]{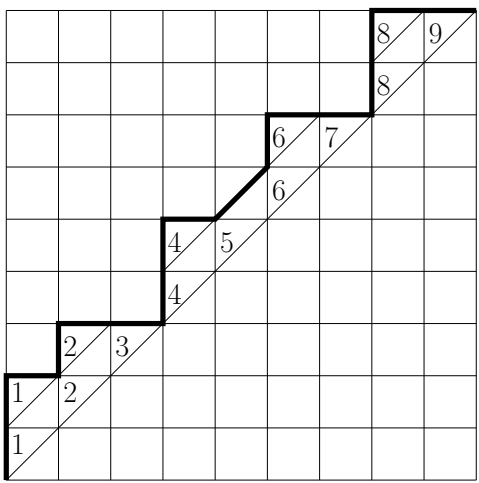}
\end{center}
\caption{\label{cap:Schroder_concat}The Schr\"oder path $p$ of Example \ref{ex:earliest_level_feature}.}
\end{figure}

\label{ex:earliest_level_feature}
\end{example}

\begin{lemma}\label{thm:latest_level_feature}
Let $p \in \D_n$ be a Schr\"{o}der path which does not contain a $\dd$ step connecting $(t-1,t-1)$ to $(t,t)$ for any $t \in \{1,2,\ldots,n\}$. If the latest level feature of $p$ is at $t=k$, then $\pi(k+1) = 1$.
\end{lemma}

\begin{proof}
Note that $\rev(p) \in \D_n$ does not contain a $\dd$ step connecting $(t-1,t-1)$ to $(t,t)$ for any $t \in \{1,2,\ldots,n\}$, and the earliest level feature of $\rev(p)$ is at $t=n-k$. Applying Lemma \ref{thm:earliest_level_feature}, and noting that $\varphi(\rev(p)) = \pi^{-1}$ by Theorem \ref{lemma:inverse_through_EM_bijection}, we obtain $\pi^{-1}(1) = k+1$, i.e., $\pi(k+1) = 1$.
\end{proof}

\begin{figure}
\begin{center}\includegraphics[%
  width=0.5\columnwidth,
  keepaspectratio]{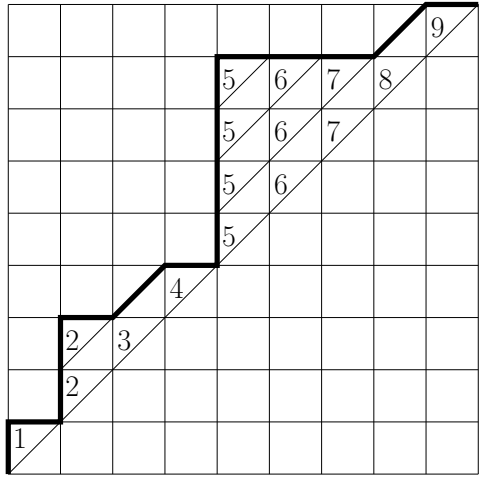}
\end{center}
\caption{\label{cap:Schroder2}The Schr\"oder path $p'$ of Example \ref{ex:example_contd}.}
\end{figure}

\begin{table}
\begin{center}
\caption{\label{ttable_concat}
Illustration of the steps involved in obtaining $\pi=\varphi(p)$ for the Schr\"oder path $p \in \mathcal{S}_9$ given in Example \ref{ex:earliest_level_feature}. The table shows the evolution of the permutation from the starting point of $(10,9,8,7,6,5,4,3,2,1)$.}

{\small
 $\begin{array}{|c|c|}
 \hline
\begin{array}[t]{c}
\mathrm{Permutation} \\
\hline 
\mathrm{Start} \\
\mathrm{permutations \; prior \; to \; that \; including} \; s_5 \\
\mathrm{after \; permutation \; including} \; s_5 \\
\mathrm{permutations \; containing} \; s_t, \; t>5 \\
\mathrm{remainder \; of \; permutations} \\
\end{array} &
\begin{array}[t]{c}
\mathrm{Result} \\
\hline 
10 \; 9 \; 8 \; 7 \; 6 \; \underline{5} \; 4 \; 3 \; 2 \; 1 \\
10 \; 9 \; 8 \; 7 \; 6 \; \underline{5} \; 4 \; \1a \; \2a \; \3a \\
10 \; 9 \; 8 \; \6a \; \underline{\5a} \; \4a \; \1a \; \7a \; 2 \; 3 \\
10 \; 9 \; 8 \; 6 \; \underline{5} \; \1a \; \4a \; 7 \; 2 \; 3 \\
\underline{\5a} \; \9a \; \8a \; \0a \; \6a \; 1 \; 4 \; 7 \; 2 \; 3 \\
\end{array} \\
\hline
\end{array}$
}
\end{center}
\end{table}
The following example illustrates Lemmas \ref{thm:earliest_level_feature} and \ref{thm:latest_level_feature} in the context of the Schr\"oder path $p \in \mathcal{S}_9$ of Example \ref{ex:running_example}.
\begin{example}
Consider the path $p \in \mathcal{S}_9$ defined in Example \ref{ex:running_example}. Note that $p \in \D_9$, and that $p$ does not contain a $\dd$ step connecting $(t-1,t-1)$ to $(t,t)$ for any $t \in \{1,2,\ldots,n\}$. The earliest and latest level features of $p$ are at $t=5$ and $t=8$ respectively. We also have $\pi(1)=5$ and $\pi(9)=1$ in accordance with Lemmas \ref{thm:earliest_level_feature} and \ref{thm:latest_level_feature} respectively.
\end{example}
\begin{theorem}\label{thm:main_transformation}
For $n \ge 0$, let $\pi \in \s_{n+1}(1243,2143)$. Then $\pi' = \rc(\pi)$ lies in $\s_{n+1}(1243,2143)$ if and only if $p = \varphi^{-1}(\pi) \in \D_n$. Furthermore, if $p \in \D_n$, then $\varphi^{-1}(\pi') = p'$ where $p' = \rev (\psi(p))$.
\end{theorem}
Before giving a proof, we provide an example in order to illustrate this Theorem.
\begin{example}\label{ex:example_contd}
Consider the path $p \in \mathcal{S}_9$ defined in Example \ref{ex:running_example}. Note that $p \in \D_9$. The path $p' = \rev(\psi(p)) = \nn\ee\nn\nn\ee\dd\ee\nn\nn\nn\nn\ee\ee\ee\dd\ee$ is illustrated in Figure \ref{cap:Schroder2}.

Here we have
$\sigma_1=s_9s_8s_7s_6s_5$, $\sigma_2=s_7s_6s_5$, $\sigma_3= s_6s_5$, $\sigma_4=s_5$, $\sigma_5=s_4s_3s_2$, $\sigma_6=s_2$ and $\sigma_7=s_1$. So
\begin{eqnarray*}
\varphi(p') &=& \sigma_7 \sigma_6 \sigma_5 \sigma_4 \sigma_3 \sigma_2 \sigma_1 (10,9,8,7,6,5,4,3,2,1) \\
&=& s_1 \; s_2 \; s_4s_3s_2 \; s_5 \; s_6s_5 \; s_7s_6s_5 \; s_9s_8s_7s_6s_5 (10,9,8,7,6,5,4,3,2,1) \\
&=& (7,10,8,2,9,3,4,5,1,6) \in \s_{10}(1243,2143)
\end{eqnarray*}
The evolution of $(10,9,8,7,6,5,4,3,2,1)$ towards $\varphi(p')$ under the application of the permutations $\sigma_i$ for $i=1,2,\ldots,k$ is shown in Table \ref{ttable_2}. The reader may check that $\varphi(p') = \rc(\pi)$, as expected.
\end{example}

\begin{lemma}
If $p$ is any Schr\"oder path of length $n \ge 0$ and $\pi = \varphi(p)$, then $\rc(\pi) \in \s_{n+1}(1243,2143)$ implies that $p \in \D_n$.
\label{lemma:one_direction}
\end{lemma}
\begin{table}
\begin{center}
\caption{\label{ttable_2}
Illustration of the bijection $\varphi(p')$ for the Schr\"oder path $p' \in \mathcal{S}_9$ given in Example \ref{ex:example_contd}. The table shows the evolution of the permutation from $(10,9,8,7,6,5,4,3,2,1)$ as each permutation $\sigma_i$ is applied for $i=1,2,\ldots,k$.}

{\small
 $\begin{array}{|c|c|}
 \hline
\begin{array}[t]{c}
\mathrm{Permutation} \\
\hline 
\mathrm{Start} \\
\sigma_1 = s_9 s_8 s_7 s_6 s_5 \\
\sigma_2 = s_7 s_6 s_5 \\
\sigma_3 = s_6 s_5 \\
\sigma_4 = s_5 \\
\sigma_5 = s_4 s_3 s_2 \\
\sigma_6 = s_2 \\
\sigma_7= s_1 \\
\end{array} &
\begin{array}[t]{c}
\mathrm{Result} \\
\hline 
10 \; 9 \; 8 \; 7 \; 6 \; 5 \; 4 \; 3 \; 2 \; 1 \\
10 \; 9 \; 8 \; 7 \; \5a \; \4a \; \3a \; \2a \; \1a \; \6a \\
10 \; 9 \; 8 \; 7 \; \4a \; \3a \; \2a \; \5a \; 1 \; 6 \\
10 \; 9 \; 8 \; 7 \; \3a \; \2a \; \4a \; 5 \; 1 \; 6 \\
10 \; 9 \; 8 \; 7 \; \2a \; \3a \; 4 \; 5 \; 1 \; 6 \\
10 \; \8a \; \7a \; \2a \; \9a \; 3 \; 4 \; 5 \; 1 \; 6 \\
10 \; \7a \; \8a \; 2 \; 9 \; 3 \; 4 \; 5 \; 1 \; 6 \\
\7a \; \0a \; 8 \; 2 \; 9 \; 3 \; 4 \; 5 \; 1 \; 6 \\
\end{array} \\
\hline
\end{array}$
}
\end{center}
\end{table}
\begin{proof}
Note that since $\pi \in \s_{n+1}(1243,2143)$ and $\rc(\pi) \in \s_{n+1}(1243,2143)$, this implies that $\pi$ avoids the pattern $\rc(1243) = 2134$ also. Suppose that points $(\alpha,\alpha+1)$ and $(\beta,\beta+1)$ lie on $p$ for some $0 \le \alpha < \beta \le n-1$, and that no part of $p$ lies on or below the line segment joining these two points. [In particular, this means that no point $(t,t+1)$ lies on $p$ for $\alpha < t < \beta$; note that this implies an occurrence of $\nn$ joining $(\alpha,\alpha+1)$ to $(\alpha,\alpha+2)$, and an occurrence of $\ee$ joining $(\beta-1,\beta+1)$ to $(\beta,\beta+1)$.] Then $\varphi(p)$ is obtained by applying a sequence of permutations $\sigma_i$ to $(n+1,n,n-1,\ldots,2,1)$ as per (\ref{eq:Egge_Mansour_bijection}). A subsequence of these permutations is as follows: first a permutation $P_1$ is applied which contains $s_{\beta+1} s_{\beta} \cdots s_{\alpha+1}$ as a subsequence, later the permutation $P_2 = s_{\beta} s_{\beta-1} \cdots s_{\alpha+1}$ is applied, and immediately after application of $P_2$, a permutation $P_3$ is applied which consists of a sequence of transpositions from $\{ s_{\alpha+1}, s_{\alpha+2}, \ldots, s_{\beta-1} \}$. The permutation $P_1$ moves a number $q > n-\alpha$ to the right of position $\beta+1$. Permutation $P_2$ moves the number $n-\alpha$ to position $\beta+1$. At this point, the numbers in positions $\alpha+1, \alpha+2, \ldots, \beta$ are in \emph{decreasing} order. If after application of $P_3$ any two of these numbers (say $f > g$) end up in decreasing order, the subsequence $(f, g, n-\alpha, q)$ would be an occurrence of the pattern $2134$, which is a contradiction. Therefore the numbers in positions $\alpha+1, \alpha+2, \ldots, \beta$ must end up in \emph{increasing order}, i.e., the points $(\alpha,\alpha+1)$ and $(\beta,\beta+1)$ on $p$ must be joined\footnote{Note that reversal of ordering of the numbers in positions $\{m, m+1, \ldots, n\}$ is effected through application of the sequence of permutations $\sigma_{m} \sigma_{m+1} \cdots \sigma_{n-1}$ where $\sigma_i = s_i s_{i-1} \cdots s_m$ for each $i=m,m+1,\ldots,n-1$; this corresponds to the path $\nn^{n-m} \ee^{n-m}$.} by an occurrence of $\nn^{\beta-\alpha} \ee^{\beta-\alpha}$. We conclude that for every pair of points $(\alpha,\alpha+1)$ and $(\beta,\beta+1)$ on $p$ such that no part of $p$ lies on or below the line segment joining these two points, the two points must be joined by an occurrence of $\nn^{\beta-\alpha} \ee^{\beta-\alpha}$. This is equivalent to the condition $p \in \D_n$.
\end{proof}
The following example illustrates Lemma \ref{lemma:one_direction} in the context of the Schr\"oder path $p \in \mathcal{S}_9$ of Example \ref{ex:running_example}.
\begin{example}
Consider the Schr\"oder path $p \in \mathcal{S}_9$ of Example \ref{ex:running_example}. In particular, note that points $(1,2)$ and $(4,5)$ lie on $p$, while points $(2,3)$ and $(3,4)$ do not; therefore we may apply the reasoning in the proof of Lemma \ref{lemma:one_direction} with $\alpha=1$ and $\beta=4$. After application of permutations $\sigma_1$ and $\sigma_2$, the permutation $(10,9,8,7,6,5,4,3,2,1)$ has changed to $(10,9,8,7,6,5,2,3,1,4)$ (c.f. row 3 of Table \ref{ttable_1}). Here $P_1=\sigma_3=s_6 s_5 s_4 s_3 s_2$ contains $s_{\beta+1} s_{\beta} \cdots s_{\alpha+1} = s_5 s_4 s_3 s_2$ as a subsequence. This permutation moves the number $q = 9 > 8 = n-\alpha$ (via a sequence of adjacent transpositions) to position $7$ which lies to the right of position $\beta+1=5$. Then, permutation $P_2 = \sigma_4 = s_{\beta} s_{\beta-1} \cdots s_{\alpha+1} = s_4 s_3 s_2$ is applied, which moves the number $n-\alpha=8$ to position $\beta+1 = 5$. To avoid the subsequence $x y 8 9$ being an occurrence of the pattern $2134$ for some $x$ and $y$, the subsequence $7  6  5$ must next be rearranged as $5  6  7$; this requires $\sigma_5= s_3 s_2$ and $\sigma_6 = s_2$, i.e., the points $(1,2)$ and $(4,5)$ on $p$ must be joined by an occurrence of $\nn^{3} \ee^{3}$.
\label{ex:must_have_spike}
\end{example}

\begin{proof}[\bf{Proof of Theorem \ref{thm:main_transformation}}]

First suppose that the Schr\"oder path $p$ contains an occurrence of $\dd$ joining $(t-1,t-1)$ to $(t,t)$ for $t \in \{ 1,2,\ldots,n \}$, i.e., the path $p$ may be viewed as the concatenation of a Schr\"oder path $q$ (of length $t-1$), a $\dd$ step, and a Schr\"oder path $r$ (of length $n-t$). From Lemma \ref{lemma:no_mixing}, this occurs if and only if the largest $t$ numbers in $\{ 1,2,\ldots, n+1\}$ occupy the first $t$ positions of $\pi = \varphi(p)$. From the definition of the reverse-complement, this latter condition is obtained if and only if the largest $n+1-t$ numbers in $\{ 1,2,\ldots, n+1\}$ occupy the first $n+1-t$ positions of $\pi' = \rc(\pi)$. Again invoking Lemma \ref{lemma:no_mixing}, this occurs if and only if the Schr\"oder path $p' = \varphi^{-1}(\pi')$ contains an occurrence of $\dd$ joining $(n-t,n-t)$ to $(n-t+1,n-t+1)$. The problem is then seen to reduce to proving the proposition for $q$ and $r$ separately (note that $p\in \D_n$ if and only if both $q \in \D_{t-1}$ and $r \in \D_{n-t}$). For this reason, in the following we need consider only those Schr\"oder paths $p$ which do not contain any occurrence of $\dd$ joining $(t-1,t-1)$ to $(t,t)$ for $t \in \{ 1,2,\ldots,n \}$. By Lemma \ref{thm:latest_level_feature}, if the latest level feature of such a path $p$ is at $t=k$, then $\pi(k+1) = 1$.

From Lemma \ref{lemma:one_direction}, if $p$ is any Schr\"oder path of length $n \ge 0$ and $\pi = \varphi(p)$, then $\rc(\pi) \in \s_{n+1}(1243,2143)$ implies that $p \in \D_n$. To prove the other direction, let $p \in \D_n$, and let $\varphi(p) = \pi$. Let $p' = \rev (\psi(p))$ and $\varphi(p') = \pi'$. We wish to show that $\pi' = \rc(\pi)$. We proceed by induction on the number $\gamma$ of features of $p$. The base case $\gamma = 0$ of the induction has already been proved in Lemma \ref{lemma:induction_base_case}.

Next assume the result holds for all Schr\"oder paths with $\gamma \ge 0$ features, and consider a Schr\"oder path $p$ with $\gamma+1$ features. The latest feature of $p$ occurs at $t=n-i$ for some integer $i\in\{ 1,2,\ldots, n-1 \}$. We assume also that this is a notch feature (this assumption will be justified later). It follows that we may write $p$ as the concatenation of two Schr\"oder paths: $q$ (of length $n-i$, and containing $\gamma$ features) and $r$ (of length $i$, and containing no features). 

The permutation $\pi_2=\varphi(q)$ is obtained by applying some sequence of permutations $\{ \tau_j \}$ to the permutation $(n+1-i,n-i,n-1-i,\ldots, 2,1)$. Similarly, the permutation $\pi_1=\varphi(r)$  is obtained by applying some sequence of permutations $\{ \sigma_j \}$ to the permutation $(i+1,i,i-1,\ldots, 2,1)$. The image of $p$ under the bijection $\varphi$ is then obtained by applying the sequence of permutations $\{ \nu_j \}$, followed by the sequence of permutations $\{ \tau_j \}$, to the permutation $(n+1,n,n-1,\ldots, 2,1)$, where the sequence $\{ \nu_j \}$ is simply the sequence $\{ \sigma_j \}$ with each transposition $s_t$ replaced by $s_{t+n-i}$. 

The application of the permutations $\{ \nu_j\}$ (corresponding to the path $r$) results in the final $i+1$ values of $(n+1,n,n-1,\ldots, 2,1)$ being replaced with the permutation $\pi_1$. Suppose that the latest level feature of $q$ is at $t=k$; by Lemma \ref{thm:latest_level_feature}, this implies that $\pi_2(k+1) = 1$. Then, since $\pi_1(1)=1$ (by Lemma \ref{lemma:induction_base_case}), the application of the permutations $\{ \tau_j \}$ (corresponding to the path $q$) results in the initial $n+1-i$ values of the resulting permutation being replaced with $\pi_2+i$, except for $\pi(k+1)$ which takes the value $1$.

Formally, we may summarize these results as
\begin{equation}
\varphi(p) = (f(1), f(2), \ldots, f(n+1-i), \pi_1(2), \pi_1(3), \ldots, \pi_1(i+1))
\label{eq:main_result_proof_1}
\end{equation}
where, for each $t \in \{ 1,2, \ldots, n+1-i\}$,
\begin{equation}
f(t) = \left\{ \begin{array}{cc}
1 & \textrm{ if } t = k+1 \\
\pi_2(t) + i & \textrm{ otherwise. }\end{array}\right. \;
\label{eq:ft_definition}
\end{equation}

Next consider $p' = \rev (\psi(p))$; this Schr\"oder path may be considered as the concatenation of Schr\"oder paths $r' = \rev (\psi(r))$ and $q' = \rev (\psi(q))$, followed by the replacement of the resulting notch feature at $t=i$ by a level feature. Let $\pi_2' = \varphi(q')$ and $\pi_1' = \varphi(r')$.

The permutation $\pi_2'=\varphi(q')$ is obtained by applying some sequence of permutations $\{ \tau_j' \}$ to the permutation $(n+1-i,n-i,n-1-i,\ldots, 2,1)$, and the permutation $\pi_1'=\varphi(r')$  is obtained by applying some sequence of permutations $\{ \sigma_j' \}$ to the permutation $(i+1,i,i-1,\ldots, 2,1)$. The permutation $\pi=\varphi(p')$ is then obtained by applying a sequence of permutations $\{ \mu_j' \}$ to the permutation $(n+1,n,n-1,\ldots, 2,1)$, followed by the sequence of permutations $\{ \sigma_j' \}$ with $\sigma'_1 = s_i s_{i-1} \cdots s_2 s_1$ omitted. The application of the permutations in $\{ \mu_j' \}$ results in the final $n+1-i$ values of the resulting permutation being replaced with $\pi_2'$, except for position $n+1-k$ which holds the value $n+1$. The reason for the latter condition is that the level feature at $t=i$ causes the number $n+1$ to move into position $n+1-k$. The reason it moves into this particular position is that the latest level feature of $q$ lying at $t=k$ implies that the earliest notch feature of 
$p'$ is at $t = n-k$. 

Note that at this point, the first $i$ values are in decreasing order starting with $n$. Therefore, application of the permutations $\{ \sigma_j'\}$ with $\sigma'_1 = s_i s_{i-1} \cdots s_2 s_1$ omitted (corresponding to the path $r'$) results in the first $i$ values of the resulting permutation being replaced with the first $i$ numbers in $\pi_1'+n-i$. 

Formally, we may summarize these results as
\begin{align}
\varphi(p') &= (\pi'(1) + n-i,\pi'(2) + n-i, \ldots, \pi'(i) + n-i, \nonumber \\
& \quad \quad \quad \quad \quad \quad \quad \quad g(1), g(2), \ldots, g(n+1-i)) \; ,
\label{eq:main_result_proof_2}
\end{align}
where, for each $t \in \{ 1,2, \ldots, n+1-i\}$,
\begin{equation}
g(t) = \left\{ \begin{array}{cc}
n+1 & \textrm{ if } t = n+1-i-k \\
\pi_2'(t) & \textrm{ otherwise. }\end{array}\right. \;
\label{eq:gt_definition}
\end{equation}
Comparing (\ref{eq:main_result_proof_1}) and (\ref{eq:ft_definition}) with (\ref{eq:main_result_proof_2}) and (\ref{eq:gt_definition}) while applying the induction hypothesis, which guarantees $\pi_1' = \rc(\pi_1)$ and $\pi_2' = \rc(\pi_2)$, yields $\varphi(p') = \rc(\pi)$, as required.

Finally, note that the preceding proof assumed that the latest feature of $p$ is a notch feature. If this is not the case, we may replace $p$ by $\psi(p)$ (whose latest feature \emph{is} a notch feature) and repeat the argument from the beginning, thus establishing (due to Theorem \ref{lemma:inverse_through_EM_bijection}) that $(\pi')^{-1} = \rc(\pi^{-1})$. By Lemma \ref{lemma:rc_inv}, this implies that $\pi' = \rc(\pi)$, as required. The result then follows by the principle of induction.

\end{proof}

The following example illustrates the inductive proof of Theorem \ref{thm:main_transformation} for the case of $\gamma=2$.
\begin{example}
Consider the Schr\"{o}der path $p = \nn\nn\ee\dd\ee\nn\nn\nn\ee\ee\ee\nn\nn\ee\nn\ee\ee \in \mathcal{D}_9$ with $\gamma+1=3$ features shown in Figure \ref{cap:main_result} (a). This path is a concatenation of a path $q \in \mathcal{D}_6$ and a featureless path $r \in \mathcal{D}_3$, i.e., $i=3$. The reader may verify that $\pi_1 = \varphi(r) = (1,3,2,4)$ and $\pi_2 = \varphi(q) = (5,6,1,7,2,3,4)$. Here the latest level feature of $q$ is at $k=2$, so Lemma \ref{thm:latest_level_feature} guarantees $\pi_2(k+1) = \pi_2(3) = 1$. The application of the first set of permutations $\{ \nu_j \}$ (corresponding to the path $r$) results in the final $i+1=4$ values of $(10,9,8,7,6,5,4,3,2,1)$ being replaced with $\pi_1$ (c.f. row 2 of upper Table \ref{ttable_4}). Then, since $\pi_1(1)=1$, the application of the second set of permutations $\{ \tau_j \}$ (corresponding to the path $q$) results in the initial $n+1-i=7$ values of the resulting permutation being replaced with $\pi_2+i = \pi_2+3 = (8,9,4,10,5,6,7)$, except for $\pi(k+1) = \pi(3)$ which takes the value $1$ (c.f. row 3 of upper Table \ref{ttable_4}).

\begin{figure}
\begin{center}
  \includegraphics[%
  width=0.48\columnwidth,
  keepaspectratio]{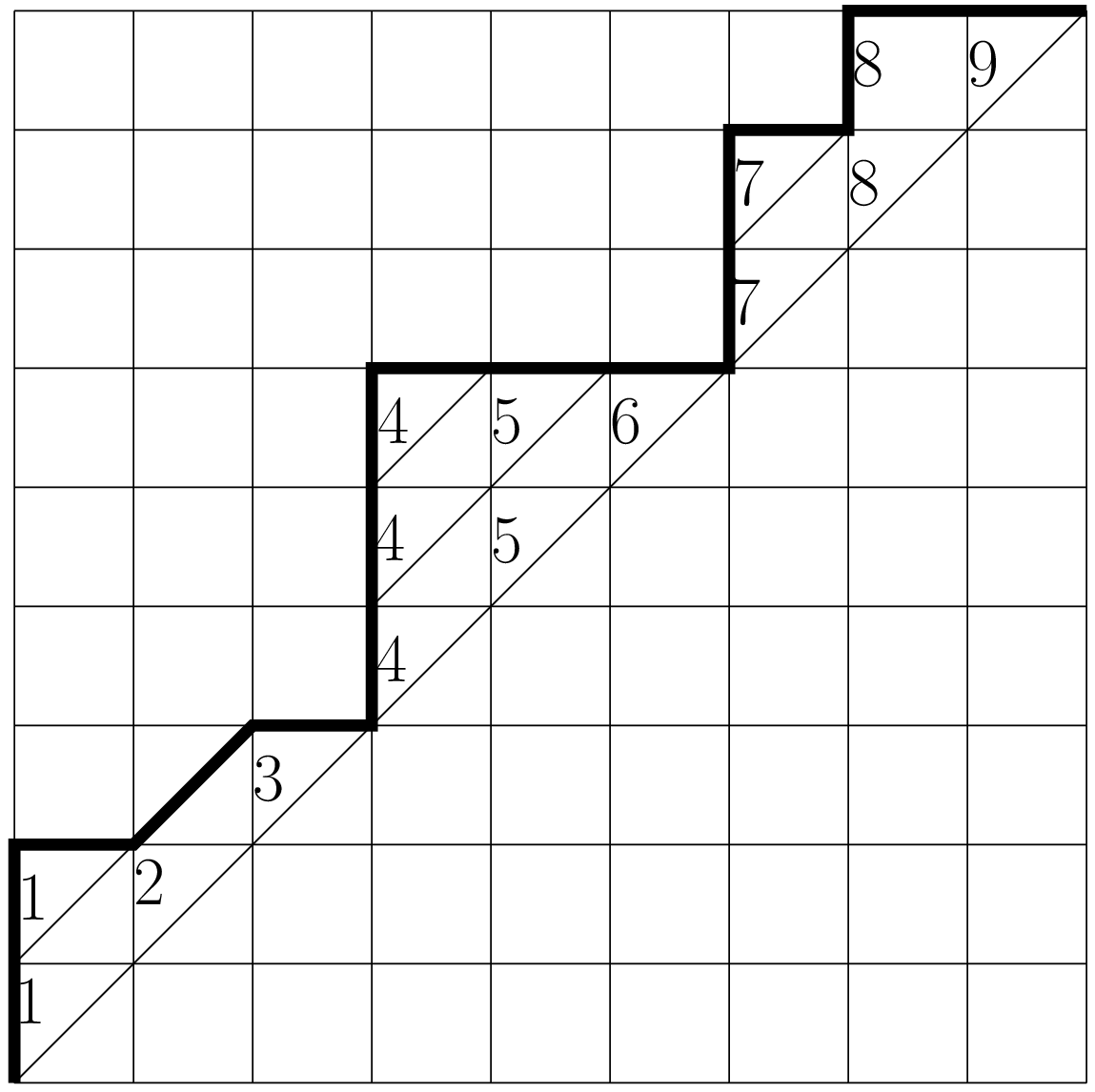}
  \includegraphics[%
  width=0.48\columnwidth,
  keepaspectratio]{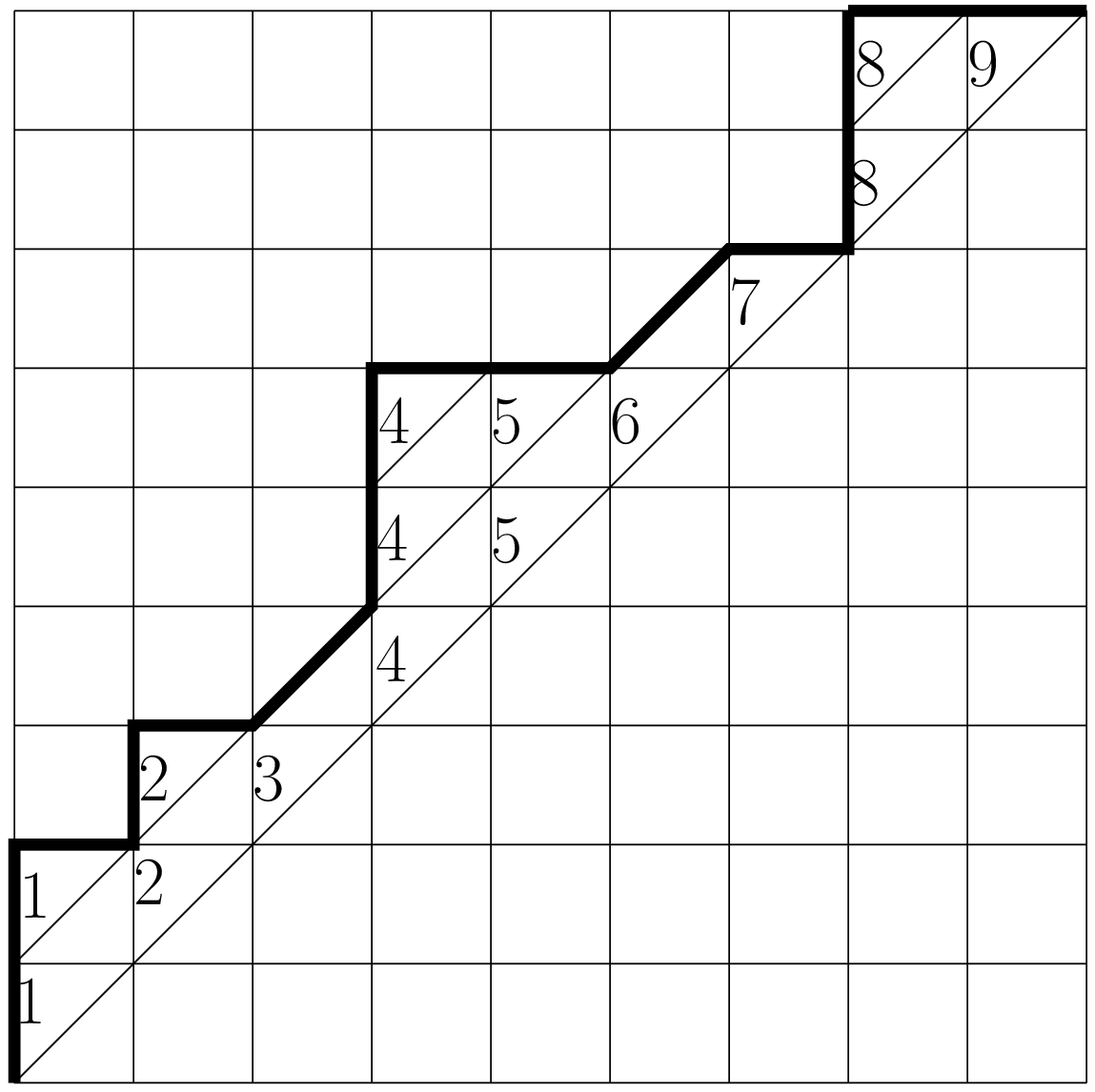}
\vspace{2mm} \\
(a) \hspace{5cm} (b) 
\end{center}
\caption{\label{cap:main_result} (a) The Schr\"oder path $p$ and (b) the Schr\"oder path $p' = \rev (\psi(p))$ of Example \ref{ex:main_result}. The corresponding evolution of permutations via the bijection $\varphi$ is shown in Table \ref{ttable_4}.}
\end{figure}

The path $p' = \rev(\psi(p))$ is shown in Figure \ref{cap:main_result} (b); this may be seen as the concatenation of $r'=\rev(\psi(r))$ and $q'=\rev(\psi(q))$, followed by the replacement of the resulting notch feature at $i=3$ by a level feature. The reader may verify that $\pi_1' = \psi(r') = (1,3,2,4)$ and $\pi_2' = \psi(q') = (4,5,6,1,7,2,3)$. The application of the first set of permutations $\{ \mu_j' \}$ results in (i) the number $n+1 = 10$ moving into position $n+1-k = 8$ (c.f. row 2 of lower Table \ref{ttable_4}) -- the reason it moves into this position is that the latest level feature of $q$ lying at $k=2$ implies that the earliest notch feature of $p'$ lies at $n-k=7$; and (ii) the final $n+1-i=7$ values of the resulting permutation being replaced with $\pi_2'$, except for position $n+1-k = 8$ which holds the value $n+1=10$. At this point, the first $i=3$ values are in decreasing order starting with $n=9$ (again c.f. row 2 of lower Table \ref{ttable_4}). Application of the permutations $\{ \sigma_j'\}$ with $\sigma_1 = s_3 s_2 s_1$ omitted (corresponding to the path $r'$) results in the first $i=3$ values of the resulting permutation being replaced with the first $i=3$ numbers in $\pi_1'+n-i = \pi_1'+6 = (7,9,8,10)$ (c.f. row 3 of lower Table \ref{ttable_4}). Application of the induction hypothesis, which guarantees $\pi_1' = \rc(\pi_1)$ and $\pi_2' = \rc(\pi_2)$, then yields the result. 
\label{ex:main_result}
\end{example}

\begin{theorem}\label{Schroder_perms}
\[
\left| \C_{2n}(1243,2143) \right| = \left| \C_{2n+1}(1243,2143) \right| = q_n
\]
for all $n \ge 1$, where the sequence $q_n$ is defined by $q_1=2$, $q_2=7$ and for every $n \ge 3$, $q_n = 4q_{n-1} - q_{n-2}$.
\end{theorem}

\begin{table}
\begin{center}
\caption{\label{ttable_4}
Illustration of the steps involved in obtaining $\pi=\varphi(p)$ (upper) and $\pi' = \varphi(p')$ (lower) for the Schr\"oder paths $p,p' \in \mathcal{S}_9$ given in Example \ref{ex:main_result}. The table shows the evolution of the permutation from $(10,9,8,7,6,5,4,3,2,1)$ in each case.}

{\small
 $\begin{array}{|c|c|}
 \hline
\begin{array}[t]{c}
\mathrm{Permutation} \\
\hline 
\mathrm{Start} \\
\{ \nu_j \} \\
\{ \tau_j \} \\
\end{array} &
\begin{array}[t]{c}
\mathrm{Result} \\
\hline 
10 \; 9 \; 8 \; 7 \; 6 \; 5 \; 4 \; 3 \; 2 \; 1 \\
10 \; 9 \; 8 \; 7 \; 6 \; 5 \; \underline{\1a} \; \3a \; \2a \; \4a \\
\8a \; \9a \; \underline{\1a} \; \0a \; \5a \; \6a \; \7a \; 3 \; 2 \; 4 \\
\end{array} \\
\hline
\end{array}$
}
{\small
 $\begin{array}{|c|c|}
 \hline
\begin{array}[t]{c}
\mathrm{Permutation} \\
\hline 
\mathrm{Start} \\
\{ \mu_j' \} \\
\{ \sigma'_j \} \\
\end{array} &
\begin{array}[t]{c}
\mathrm{Result} \\
\hline 
10 \; 9 \; 8 \; 7 \; 6 \; 5 \; 4 \; 3 \; 2 \; 1 \\
\9a \; \8a \; \7a \; \4a \; \5a \; \6a \; \1a \; \underline{\0a} \; \2a \; \3a \\
\7a \; \9a \; \8a \; 4 \; 5 \; 6 \; 1 \; \underline{10} \; 2 \; 3 \\
\end{array} \\
\hline
\end{array}$
}
\end{center}
\end{table}

\begin{proof}
Denote by $\D$ the set of all Schr\"{o}der prefixes which may be extended to form Schr\"{o}der paths in $\D_n$ for any finite $n$. For $i \ge 0$, let $a_i$ denote the number of Schr\"oder prefixes in $\mathcal{D}$ terminating at the point $(i,i)$ (this is simply $|D_i|$). For $i \ge 1$, let $b_i$ denote the number of Schr\"oder prefixes in $\mathcal{D}$ terminating at the point $(i-1,i)$. We have $a_0 = 1$; also for completeness we define $b_0 = 0$. For $i \ge 1$ the point $(i,i)$ may be reached either by a $\dd$ step or by an $\ee$ step, and so we have 
\begin{equation}
a_i = a_{i-1} + b_i \quad \mbox{for} \: i \ge 1 \; .
\label{eq:a_recursion}
\end{equation} 
The justification of this recursion is illustrated in Figure \ref{cap:paths_and_recursion2} for the case of $i=7$. Summing (\ref{eq:a_recursion}) over $i=1,2,\ldots, k$ we obtain 
\begin{equation}
a_k = 1 + \sum_{i=0}^{k} b_i \quad \mbox{for} \: k \ge 0 \; .
\label{eq:a_recursion2}
\end{equation}
We have $b_1 = a_0$ since the point $(0,1)$ may only be reached by an $\nn$ step. For $i \ge 1$ consider the point $(i-1,i)$; there are $a_{i-1}$ paths which reach this point via an $\nn$ step, $b_{i-1}$ paths which reach this point via a $\dd$ step, and $b_j$ paths which reach this point via the steps ${\nn}^{i-j}{\ee}^{i-j}$ for each $j=1,2,\ldots, i-1$. Therefore, for each $i \ge 1$,
\begin{equation}
b_{i} = a_{i-1} + b_{i-1} + \sum_{j=1}^{i-1} b_j = 1 + b_{i-1} + 2\sum_{j=0}^{i-1} b_j
\label{eq:recursion_bis}
\end{equation}
where the second equality is obtained using (\ref{eq:a_recursion2}). The justification of this recursion is illustrated in Figure \ref{cap:paths_and_recursion2} for the case of $i=5$. Subtracting (\ref{eq:recursion_bis}) for $i=k$ from (\ref{eq:recursion_bis}) for $i=k+1$ we obtain
\begin{equation}
b_{k+1} = 4b_k - b_{k-1}
\label{eq:b_recursion}
\end{equation}
for all $k \ge 1$, with $b_0=0$ and $b_1 = 1$.

\begin{figure}
\begin{center}
  \includegraphics[%
  width=0.52\columnwidth,
  keepaspectratio]{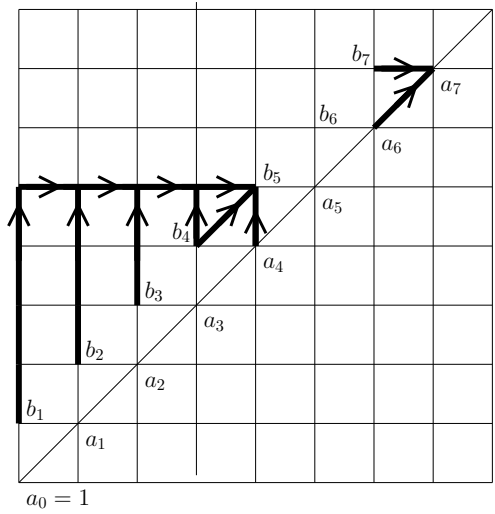}
\end{center}
\caption{\label{cap:paths_and_recursion2} The figure shows graphically the justification of (\ref{eq:a_recursion}) and (\ref{eq:recursion_bis}). In the figure, each lattice point $(i,i)$ (resp. $(i,i-1)$) is labeled with the number $a_i$ (resp. $b_i$) of Schr\"{o}der prefixes in $\mathcal{D}$ which terminate at that lattice point. The figure shows all possible terminations of such prefixes at the points $(7,7)$ and $(4,5)$, thus illustrating that (respectively) $a_7 = a_6 + b_7$ and $b_5 = a_4 + 2 b_4 + b_3 + b_2 + b_1$.}
\end{figure}
\begin{figure}
\begin{center}
  \includegraphics[%
  width=0.75\columnwidth,
  keepaspectratio]{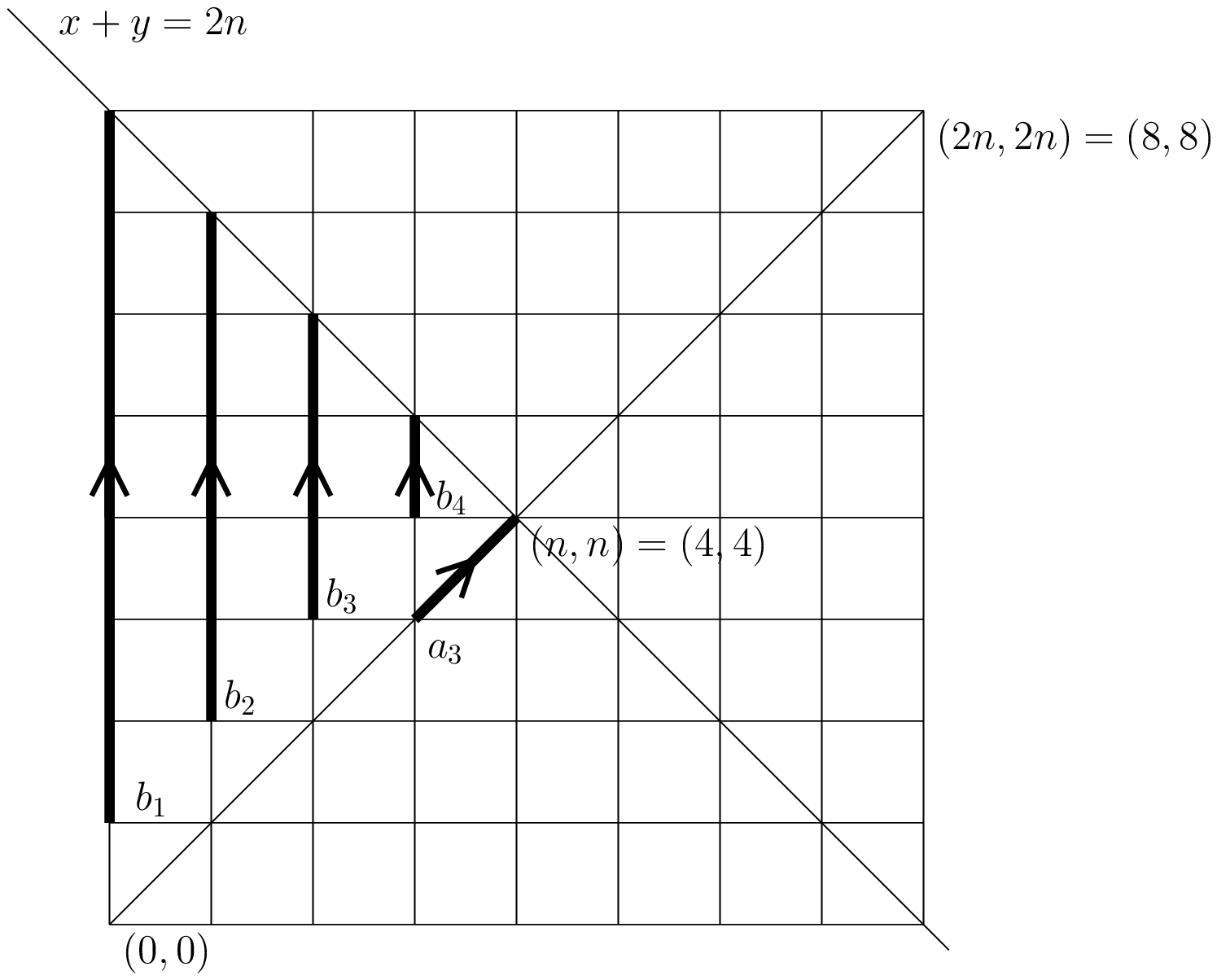}
\end{center}
\caption{\label{cap:paths_and_recursion1} The figure shows graphically the justification of (\ref{eq:qn_recursion}) for the example of $n=4$. Any Schr\"{o}der prefix in $\mathcal{D}$ which terminates on the line $x+y=2n$ must have one of the possible terminations shown in the figure. Each such prefix then has a unique completion to form a Schr\"{o}der path $p \in \D_{2n}$. Therefore, in the case illustrated, $q_4 = a_3 + b_4 + b_3 + b_2 + b_1$. A similar illustration may be made for the case of odd Schr\"{o}der path length.}
\end{figure}

Next, by Theorem \ref{thm:main_transformation}, the number of permutations $\pi \in \s_{n+1}(1243,2143)$ which satisfy $\rc(\pi) = \pi$ is equal to the number of Schr\"oder paths $p \in \mathcal{D}_{n}$ which satisfy $\rev(\psi(p)) = p$. First let $q_n$ denote the number of Schr\"oder paths $p \in \mathcal{D}_{2n}$ which satisfy $\rev(\psi(p)) = p$. The initial steps of any such path must form a Schr\"oder prefix in $\mathcal{D}$ terminating on the line $x+y=2n$. This termination occurs either at the point $(n,n)$ (there are $a_{n-1}$ of these -- note that the point $(n,n)$ may not be reached by an $\ee$ step) or at the point $(i,2n-i)$ for some $i \in \{0,1,\ldots, n-1 \}$ (there are $b_{i+1}$ of these, as they must join $(i,i+1)$ to $(i,2n-i)$ via the steps $\nn^{2(n-i)-1}$). Note that the point $(n,n+1)$ may not be reached via a $\dd$ step. It is easy to see that each of these Schr\"oder prefixes has a unique completion to form a Schr\"oder path with $\rev(\psi(p)) = p$. Therefore 
\begin{equation}
q_n = a_{n-1} + \sum_{i=0}^{n} b_i = 2b_n - b_{n-1}
\label{eq:qn_recursion}
\end{equation}
for $n \ge 1$ (using (\ref{eq:a_recursion2}) and (\ref{eq:recursion_bis})). The justification of this recursion is illustrated in Figure \ref{cap:paths_and_recursion1}. From (\ref{eq:b_recursion}) we then have $q_{n+1} = 4q_n - q_{n-1}$ with $q_1=2$ and $q_2=7$.

Similarly, let $u_n$ denote the number of Schr\"oder paths $p \in \mathcal{S}_{2n-1}$ which satisfy $\rev(\psi(p)) = p$. The initial steps of any such path must form a Schr\"oder prefix terminating at the point $(n-1,n-1)$ (with the next step joining $(n-1,n-1)$ to $(n,n)$ via a $\dd$ step -- there are $a_{n-1}$ of these), or at the point $(i,2n-1-i)$ for some $i \in \{0,1,\ldots, n-1 \}$ (there are $b_{i+1}$ of these, as they must join $(i,i+1)$ to $(i,2n-1-i)$ via the steps $\nn^{2(n-i-1)}$). Again, each of these Schr\"oder prefixes has a unique completion to form a Schr\"oder path with $\rev(\psi(p)) = p$; thus $u_n = a_{n-1} + \sum_{i=0}^{n} b_i = 2b_n - b_{n-1}$ for $n \ge 1$ and so $u_n=q_n$ for $n \ge 1$.
\end{proof}

\begin{theorem}\label{Schroder_invs}
\[
\left| \CI_{2n}(1243,2143) \right| = \left| \CI_{2n+1}(1243,2143) \right| = p_n
\]
for all $n \ge 1$, where $p_n$ denotes the $n$-th Pell number, i.e., $p_1=2$, $p_2=5$ and for every $n \ge 3$, $p_n = 2p_{n-1} + p_{n-2}$.
\end{theorem}

\begin{proof}
This proceeds similarly to the proof of Theorem \ref{Schroder_perms}. From Corollary \ref{cor:involutions}, a permutation $\pi \in \s_n(1243,2143)$ is an involution if and only if $p = \varphi^{-1}(\pi)$ is \emph{symmetric}, i.e., if and only if it satisfies $\rev(p) = p$. Therefore, the number of involutions in $\s_{n+1}(1243,2143)$ which satisfy $\rc(\pi) = \pi$ is equal to the number of Schr\"oder paths in $p \in \mathcal{S}_n$ which satisfy $\psi(p) = p$; our task is to count these Schr\"oder paths. To this end, let $\tilde{\mathcal{D}}$ denote the set of Schr\"oder prefixes in $\mathcal{D}$ with no features.
For $i \ge 0$, let $c_i$ denote the number of Schr\"{o}der prefixes in $\tilde{\mathcal{D}}$ terminating at the point $(i,i)$. For $i \ge 1$, let $d_i$ denote the number of prefixes in $\tilde{\mathcal{D}}$ terminating at the point $(i-1,i)$. We have $c_0 = 1$; also for completeness define $d_0 = 0$. For $i \ge 1$ the point $(i,i)$ may be reached either by a $\dd$ step or by an $\ee$ step, and so we obtain 
\begin{equation}
c_k = 1 + \sum_{i=0}^{k} d_i \quad \mbox{for} \: k \ge 0
\label{eq:c_recursion}
\end{equation}
by the same method as that which obtained (\ref{eq:a_recursion2}). We have $d_1 = c_0$ since the point $(0,1)$ may only be reached by an $\nn$ step. For $i \ge 2$ consider the point $(i-1,i)$; there are $c_{i-2}$ paths which reach this point via an $\nn$ step (since such a step must be preceded by a $\dd$ step), and $d_j$ paths which reach this point via the steps ${\nn}^{i-j}{\ee}^{i-j}$ for each $j=1,2,\ldots, i-1$. Therefore, for each $i \ge 2$,
\begin{equation}
d_{i} = c_{i-2} + d_{i-1} + \sum_{j=1}^{i-1} d_j = 1 + d_{i-1} + 2\sum_{j=0}^{i-2} d_j
\label{eq:recursion_dis}
\end{equation}
where we have used (\ref{eq:c_recursion}). Subtracting (\ref{eq:recursion_dis}) for $i=k$ from (\ref{eq:recursion_dis}) for $i=k+1$ we obtain
\begin{equation}
d_{k+1} = 2d_k + d_{k-1}
\label{eq:d_recursion}
\end{equation}
for all $k \ge 1$, with $d_0=0$ and $d_1 = 1$.

Let $p_n$ denote the number of Schr\"oder paths $p \in \mathcal{S}_{2n}$ which satisfy $\psi(p) = p$. The initial steps of any such path must form a Schr\"oder prefix in $\tilde{\mathcal{D}}$ terminating at the point $(n,n)$ (there are $c_{n-1}$ of these -- note that the point $(n,n)$ may not be reached by an $\ee$ step) or at the point $(i,2n-i)$ for some $i \in \{0,1,\ldots, n-1 \}$ (there are $d_{n+1}$ of these, as they must join $(i,i+1)$ to $(i,2n-i)$ via the steps $\nn^{2(n-i)-1}$). It is easy to see that each of these Schr\"oder prefixes has a unique completion to form a Schr\"oder path with $\psi(p) = p$. Therefore $p_n = c_{n-1} + \sum_{i=0}^{n} d_i = d_{n+1}$ for $n \ge 1$ (again using (\ref{eq:c_recursion})). From (\ref{eq:d_recursion}) we then have $p_{n+1} = 4p_n - p_{n-1}$ with $p_1=2$ and $p_2=5$.

Similarly, let $v_n$ denote the number of Schr\"oder paths $p \in \mathcal{S}_{2n-1}$ which satisfy $\psi(p) = p$. The initial steps of any such path must form a Schr\"oder prefix terminating at the point $(n-1,n-1)$ (with the next step joining $(n-1,n-1)$ to $(n,n)$ via a $\dd$ step -- there are $c_{n-1}$ of these), or at the point $(i,2n-1-i)$ for some $i \in \{0,1,\ldots, n-1 \}$ (there are $d_{n+1}$ of these, as they must join $(i,i+1)$ to $(i,2n-1-i)$ via the steps $\nn^{2(n-i-1)}$). Again, each of these Schr\"oder prefixes has a unique completion to form a Schr\"oder path with $\psi(p) = p$; thus $v_n = c_{n-1} + \sum_{i=0}^{n} d_i = d_{n+1}$ for $n \ge 1$ and so $v_n=p_n$ for $n \ge 1$.
\end{proof}

\section*{Acknowledgments}
The authors would like to thank the anonymous referees for their helpful suggestions which have greatly improved the accuracy and presentation of this work. They would also like to thank Marilena Barnabei, W. M. B. Dukes and T. Mansour for helpful discussions.

\end{document}